\documentclass[aap,MSNbibl,seceqn,dvips]{arximspdf}
\usepackage{mathrsfs}

\doi{10.1214/13-AAP926} %kopijuoti is PTS
\volume{24}
\issue{1}
\pubyear{2014}
\firstpage{337}
\lastpage{377}

\makeatletter
\newcommand{\rrVert}{\Vert}
\newcommand{\rrvert}{\vert}
\newcommand{\llVert}{\Vert}
\newcommand{\llvert}{\vert}
\newcommand{\E}{\mathbb{E}}
\renewcommand{\P}{\mathbb{P}}
\renewcommand{\tilde}{\widetilde}
\renewcommand{\epsilon}{\varepsilon}

\newtheorem{prop}{Proposition}[section]
\newtheorem{thmm}[prop]{Theorem}
\newtheorem{coroll}[prop]{Corollary}
\newtheorem{lemma}[prop]{Lemma}
\newproclaim{defi}[prop]{Definition}

\newproclaim{rmk}[prop]{Remark}
\newproclaim{ex}[prop]{Example}
\newproclaim{assumption}[prop]{Assumption}

\makeatother

\begin{document}
\begin{frontmatter}

\title{Error bounds for Metropolis--Hastings algorithms applied to
perturbations of Gaussian measures in high dimensions}
\runtitle{Error bounds for Metropolis--Hastings algorithms}

\begin{aug}
\author{\fnms{Andreas} \snm{Eberle}\corref{}\ead[label=e1]{eberle@uni-bonn.de}}
\address{Institut f\"ur Angewandte Mathematik\\
Universit\"at Bonn\\
Endenicher Allee 60\\
53115 Bonn\\
Germany\\
\printead{e1}}
\affiliation{Universit\"at Bonn}
\runauthor{A. Eberle}
\end{aug}

% HISTORY:
\received{\smonth{11} \syear{2012}}
\revised{\smonth{2} \syear{2013}}

% ABSTRACT
%
\begin{abstract}
The Metropolis-adjusted Langevin algorithm (MALA) is a Metro\-polis--Hastings method for approximate
sampling from continuous distributions. We derive upper bounds for the
contraction rate in
Kantorovich--Rubinstein--Wasserstein distance of the MALA chain with
semi-implicit Euler proposals
applied to log-concave probability measures that have a density w.r.t. a Gaussian reference measure.
For sufficiently ``regular'' densities, the estimates are
dimension-independent, and they hold for
sufficiently small step sizes $h$ that do not depend on the dimension
either. In the limit
$h\downarrow0$, the bounds approach the known optimal contraction
rates for overdamped Langevin
diffusions in a convex potential.

A similar approach also applies to Metropolis--Hastings chains with
Ornstein--Uhlenbeck proposals. In
this case, the resulting estimates are still independent of the
dimension but less optimal, reflecting
the fact that MALA is a higher order approximation of the diffusion
limit than Metropolis--Hastings
with Ornstein--Uhlenbeck proposals.
\end{abstract}

% KEYWORDS



\end{frontmatter}

%s1 #&#
\section{Introduction}
The performance of Metropolis--Hastings (MH) methods \cite{MRTT,H,RC} for
sampling probability measures on high-dimensional continuous state
spaces has attracted
growing attention in recent years. The pioneering works by Roberts,
Gelman and Gilks
\cite{RGG} and Roberts and Rosenthal \cite{RR} show in particular that
for product
measures $\pi^d$ on $\mathbb R^d$, the average acceptance probabilities
for the
Random Walk Metropolis algorithm (RWM) and the Metropolis adjusted Langevin
algorithm (MALA) converge to a strictly positive limit as $d\to\infty$
only if the
step sizes $h$ go to zero of order $O(d^{-1 })$, $O(d^{-1/3})$,
respectively. In
this case, a diffusion limit as $d\to\infty$ has been derived, leading
to an optimal scaling
of the step sizes maximizing the speed of the limiting diffusion, and
an asymptotically
optimal acceptance probability.

Recently, the optimal scaling results for RWM and MALA have
been extended significantly to targets that are not of product form but
have a
sufficiently regular density w.r.t. a Gaussian measure; cf. \cite
{MPS,PSTb}.
On the other hand, it has been pointed out
\cite{BRSV,BS,HSV,CRSW} that for corresponding
perturbations of Gaussian measures, the
acceptance probability has a strictly positive limit as $d\to\infty$
for small step sizes
that do not depend on the dimension, provided the random walk or Euler proposals
in RWM and MALA are replaced by Ornstein--Uhlenbeck or semi-implicit
(``preconditioned'') Euler proposals, respectively; cf. also below.
Pillai, Stuart and
Thi\'ery \cite{PST} show that in this case, the Metropolis--Hastings algorithm
can be realized directly on an infinite-dimensional Hilbert space
arising in the limit as
$d\to\infty$, and the corresponding Markov chain converges weakly to
an infinite-dimensional overdamped Langevin diffusion as $h\downarrow0$.

Mixing properties and convergence to equilibrium of Langevin diffusions
have been
studied intensively \cite{Has,MT,B,BCG,PZ,Roy}. In particular, it is
well-known that contractivity and
exponential convergence to equilibrium in Wasserstein distance can be quantified
if the stationary distribution is strictly log-concave \cite{CL,RS};
cf. also \cite{E11}
for a recent extension to the nonlog-concave case. Because of the
diffusion limit
results, one might expect that the approximating Metropolis--Hastings
chains have
similar convergence properties. However, this heuristics may also be
wrong, since
the convergence of the Markov chains to the diffusion is known only in
a weak and
nonquantitative sense.

Although there is a huge number of results quantifying the speed of
convergence to
equilibrium for Markov chains on discrete state spaces (cf.~\cite
{P,S-C} for an overview),
there are relatively few quantitative results on Metro\-polis--Hastings
chains on $\mathbb R^d$
when $d$ is large. The most remarkable exception are the well-known works
\cite{DFK,KLS,LV,LV1,LV2} which prove an upper bound for the mixing
time that is
polynomial in the dimension for Metropolis chains with ball walk
proposals for
uniform measures on convex sets and more general log-concave measures.

Below, we develop an approach to quantify Wasserstein contractivity and
convergence to equilibrium in a dimension-independent way for the
Metropolis--Hastings chains with Ornstein--Uhlenbeck and semi-implicit Euler
proposals. Our approach applies in the strictly log-concave case (or,
more generally,
if the measure is strictly log-concave on an appropriate ball) and yields
bounds for small step sizes that are very precise. The results for
semi-implicit Euler proposals
require less restrictive assumptions than those for Ornstein--Uhlenbeck
proposals,
reflecting the fact that the corresponding Markov chain is a higher
order approximation
of the diffusion.

Our results are closely related and complementary to the recent work
\cite{HSV11}, and to the
dimension-dependent geometric ergodicity results
in \cite{B-RHV}. In particular, in \cite{HSV11}, Hairer, Stuart and Vollmer
apply related methods to establish exponential convergence to equilibrium
in Wasserstein distance for Metropolis--Hastings chains with
Ornstein--Uhlenbeck proposals
in a less quantitative way, but without assuming log-concavity. In the
context of
probability measures on function spaces, the techniques developed here
are applied in
the PhD Thesis of  Gruhlke \cite{G}.

We now recall some basic facts on Metropolis--Hastings algorithms and describe
our setup and the main results. Sections \ref{secWCMH} and \ref
{secProposal} contain
basic results on Wasserstein
contractivity of Metropolis--Hastings kernels, and contractivity of the
proposal kernels.
In Sections \ref{secRP} and \ref{secDRP}, we prove bounds quantifying
rejection
probabilities and the dependence of the rejection event on the current
state for
Ornstein--Uhlenbeck and semi-implicit Euler proposals. These bounds,
combined with an upper bound for the exit probability of the corresponding
Metropolis--Hastings chains from a given ball derived in Section~\ref{secLyapunov}
are crucial for the proof of the main results in Section~\ref{secFinal}.

%s1.1 #&#
\subsection{Metropolis--Hastings algorithms}
Let $U\dvtx \mathbb R^d\to\mathbb R$ be a lower bounded measurable function
such that
\[
\mathcal Z = \int_{\mathbb R^d}\exp\bigl(-U(x)\bigr) \,dx < \infty,
\]
and let $\mu$ denote the probability measure on $\mathbb R^d$ with
density proportional to $\exp(-U)$. We use the same letter $\mu$ for
the measure and its density, that is,
%
%
%e1.1 #&#
\begin{equation}
\label{eq1} \mu(dx) = \mu(x) \,dx = \mathcal Z^{-1}\exp\bigl(-U(x)\bigr)
\,dx.
\end{equation}
Below, we view the measure $\mu$ defined by (\ref{eq1}) as a
perturbation of the standard normal distribution $\gamma^d$ in
$\mathbb
R^d$; that is, we decompose
%
%
%e1.2 #&#
\begin{equation}
U(x) = \tfrac12 |x|^2+V(x), \qquad x\in\mathbb R^d, \label{eq19}
\end{equation}
with a measurable function $V\dvtx \mathbb R^d\to\mathbb R$, and obtain the
representation
%
%e1.3 #&#
\begin{equation}
\mu(dx) = \tilde{\mathcal Z}^{-1}\exp\bigl(-V(x)\bigr)
\gamma^d(dx) \label{eq20}
\end{equation}
with normalization constant $\tilde{\mathcal Z}=\mathcal Z/(2\pi
)^{d/2}$. Here $| \cdot|$ denotes the Euclidean norm.

Note that in $\mathbb R^d$, any probability measure with a strictly
positive density can be represented as an absolutely continuous
perturbation of $\gamma^d$ as in (\ref{eq20}). In an infinite-dimensional limit, however, the density may degenerate. Nevertheless,
also on infinite-dimensional spaces, absolutely continuous
perturbations of Gaussian measures form an important and widely used
class of models.

%
%ex1.1 #&#
\begin{ex}[(Transition path sampling)]\label{exTPS1}
We briefly describe a typical application; cf. \cite{HSV} and \cite{G}
for details. Suppose
that we are interested in sampling a trajectory of a diffusion process
in $\mathbb R^\ell$
conditioned to a given endpoint $b$ at time $t=1$. We assume that the
unconditioned
diffusion process $(Y_t, \mathbb P )$ satisfies a stochastic
differential equation of the form
%
%
%e1.4 #&#
\begin{equation}
dY_t = -\nabla H (Y_t) \,dt + dB_t,
\end{equation}
where $(B_t)$ is an $\ell$-dimensional Brownian motion, and $H\in
C^2(\mathbb R^\ell)$
is bounded from below. Then, by Girsanov's theorem and It\^o's
formula, a regular
version of the law of the conditioned process satisfying $Y_0=a$ and
$Y_1=b$ on the path
space $E=\{ y\in C([0,1],\mathbb R^\ell)\dvtx y_0=a,y_1=b\}$ is given by
%
%
%e1.5 #&#
\begin{equation}
\label{eqexdot} \mu(dy) = C^{-1}\exp\bigl(-V(y)\bigr) \gamma(dy),
\end{equation}
where $\gamma$ is the law of the Brownian bridge from $a$ to $b$,
%
%
%e1.6 #&#
\begin{equation}
\label{eqexdotdot} V (y) = \frac12\int_0^1
\phi(y_s) \,ds \qquad\mbox{with } \phi(x) = \bigl|\nabla H(x)\bigr|^2 -
\Delta H (x),
\end{equation}
and $C=\exp{(H(b)-H(a))}$; cf. \cite{Roy}.
In order to obtain finite-dimensional
approximations of the measure $\mu$ on $E$, we consider the \emph
{Wiener--L\'evy expansion}
%
%
%e1.7 #&#
\begin{equation}
\label{eqexstar} y_t = e_t+\sum
_{n=0}^\infty\sum_{k=0}^{2^n-1}
\sum_{i=1}^ \ell x_{n,k,i}
e_t^{n,k,i},\qquad t\in[0,1],
\end{equation}
of a path $y\in E$ in terms of the basis functions $e_t=(1-t)a+tb$ and
$e_t^{n,k,i}=2^{-n/2}
g(2^nt-k)e^i$ with $g(s)=\min(s,1-s)^+$. Here the coefficients
$x_{n,k,i}$, $n\ge0$,
$0\le k<2^n$, $1\le i\le\ell$, are real numbers. Recall that
truncating the series at $n=m-1$
corresponds to taking the polygonal interpolation of the path $y$
adapted to the dyadic
partition $\mathcal D_m=\{ k2^{-m}\dvtx k=0,1,\ldots, 2^m\} $ of the
interval $[0,1]$.
Now fix $m\in\mathbb N$, let $d=(2^m-1)\ell$ and let
\[
x^d = \bigl(x_{n,k,i}\dvtx 0\le n<m, 0\le k<2^n, 1\le
i\le l \bigr) \in \mathbb R^d
\]
denote the vector consisting of the first $d$ components in the basis
expansion of a path
$y\in E$. Then the image of the Brownian bridge measure $\gamma$ under
the projection
$\pi_d\dvtx E\to\mathbb R^d$ that maps $y$ to $x^d$ is the $d$-dimensional
standard normal
distribution $\gamma^d$; for example, cf.  \cite{Ste}. Therefore, a
natural finite-dimensional
approximation to the infinite-dimensional sampling problem described
above consists in
sampling from the probability measure
%
%
%e1.8 #&#
\begin{equation}
\label{eqexcirc} \mu_d(dx) = \tilde{\mathcal Z}_d^{-1}
\exp\bigl(-V_d(x)\bigr) \gamma^d(dx)
\end{equation}
on $\mathbb R^d$ where $\tilde{\mathcal Z}_d$ is a normalization
constant, and
%
%
%e1.9 #&#
\begin{equation}
\label{eqexcirccirc} V_d(x) = 2^{-m-1} \Biggl( \frac12
\phi(y_0)+\sum_{k=1}^{2^m-1}\phi
(y_{k2^{-m}})+\frac12 \phi(y_1) \Biggr);
\end{equation}
with $y=e+\sum_{n<m}\sum_k\sum_ix_{n,k,i}e^{n,k,i}$ denoting the
polygonal path corresponding to $x^d=(x_{n,k,i})\in\mathbb R^d$.
\end{ex}

Returning to our general setup, suppose that $p(x,dy)=p(x,y) \,dy$ is an
absolutely continuous transition kernel on $\mathbb R^d$ with strictly
positive densities $p(x,y)$. Let
%
%
%e1.10 #&#
\begin{equation}
\label{eq2} \alpha(x,y) = \min \biggl(\frac{\mu(y)p(y,x)}{\mu
(x)p(x,y)},1 \biggr),\qquad x,y\in\mathbb
R^d.
\end{equation}
Note that $\alpha(x,y)$ does not depend on $\mathcal Z$. The
Metropolis--Hastings algorithm with proposal kernel $p$ is the
following Markov chain Monte Carlo method for approximate sampling and
Monte Carlo integration w.r.t. $\mu$:

\begin{longlist}[(1)]
\item[(1)] Choose an initial state $X_0$.
\item[(2)] For $n:=0,1,2,\ldots$:
\begin{itemize}
\item Sample $Y_n\sim p(X_n,dy)$ and $U_n\sim\operatorname{Unif}(0,1)$
independently.
\item If $U_n<\alpha(X_n,Y_n)$, then accept the proposal, and set
$X_{n+1}:=Y_n$,\ else reject the proposal and set $X_{n+1}:=X_n$.
\end{itemize}
\end{longlist}
The algorithm generates a time-homogeneous Markov chain
$(X_n)_{n=0,1,2,\ldots}$ with initial state $X_0$ and transition kernel
%
%
%e1.11 #&#
\begin{equation}
\label{eq3} q(x,dy) = \alpha(x,y)p(x,y)\,dy+r(x)\cdot\delta_x(dy).
\end{equation}
Here
%
%
%e1.12 #&#
\begin{equation}
\label{eq4} r(x) = 1-q\bigl(x,\mathbb R^d\setminus\{x\}\bigr) = 1-
\int_{\mathbb R^d}\alpha (x,y)p(x,y)\,dy
\end{equation}
is the average rejection probability for the proposal when the Markov
chain is at~$x$. Note that $q(x,dy)$ restricted to $\mathbb
R^d\setminus\{x\}$ is again absolutely continuous with density
\[
q(x,y) = \alpha(x,y)p(x,y).
\]
Since
\[
\mu(x)q(x,y) = \alpha(x,y)\mu(x)p(x,y) = \min\bigl(\mu(y)p(y,x),\mu(x)p(x,y)
\bigr)
\]
is a symmetric function in $x$ and $y$, the kernel $q(x,dy)$ satisfies
the \emph{detailed balance condition}
%
%
%e1.13 #&#
\begin{equation}
\label{eq5} \mu(dx)q(x,dy) = \mu(dy)q(y,dx).
\end{equation}
In particular, $\mu$ is a stationary distribution for the
Metropolis--Hastings chain, and the chain with initial distribution
$\mu
$ is reversible. Therefore, under appropriate ergodicity assumptions,
the distribution of $X_n$ will converge to $\mu$ as $n\to\infty
$.

To analyze Metropolis--Hastings algorithms it is convenient to
introduce the function
%
%
%e1.14 #&#
\begin{equation}
\label{eq6} G(x,y) = \log\frac{\mu(x)p(x,y)}{\mu(y)p(y,x)} = U(y)-U(x)+\log \frac
{p(x,y)}{p(y,x)}.
\end{equation}
For any $x,y\in\mathbb R^d$,
%
%
%e1.15 #&#
\begin{equation}
\label{eq7} \alpha(x,y) = \exp\bigl(-G(x,y)^+\bigr).
\end{equation}
In particular, for any $x,y,\tilde x,\tilde y\in\mathbb R^d$,
%
%
%e1.16 #&#
%e1.17 #&#
%e1.18 #&#
\begin{eqnarray}
1-\alpha(x,y)&\leq&G(x,y)^+, \label{eq8}
\\
\bigl(\alpha(x,y)-\alpha(\tilde x,\tilde y)\bigr)^+&\leq&\bigl(G(x,y)-G(\tilde
x,\tilde y)\bigr)^-\quad\mbox{and}\label{eq9}
\\
\bigl(\alpha(x,y)-\alpha(\tilde x,\tilde y)\bigr)^-&\leq&\bigl(G(x,y)-G(\tilde
x,\tilde y)\bigr)^+. \label{eq9a}
\end{eqnarray}
The function $G(x,y)$ defined by (\ref{eq6}) can also be represented in
terms of $V$: Indeed, since
\[
\log\frac{\gamma^d(x)}{\gamma^d(y)} = \frac12\bigl(|y|^2-|x|^2\bigr),
\]
we have
%
%
%e1.19 #&#
\begin{equation}
G(x,y) = V(y)-V(x)+\log\frac{\gamma^d(x)p(x,y)}{\gamma^d(y)p(y,x)}, \label{eq21}
\end{equation}
where $\gamma^d(x)=(2\pi)^{-d/2}\exp(-|x|^2/2)$ denotes the standard
normal density in $\mathbb R^d$.

%s1.2 #&#
\subsection{Metropolis--Hastings algorithms with Gaussian proposals}
We aim at proving contractivity of Metropolis--Hastings kernels w.r.t. appropriate Kantorovich--Rubinstein--Wasserstein distances. For this
purpose, we are looking for a proposal kernel that has adequate
contractivity properties and sufficiently small rejection
probabilities. The rejection probability is small if the proposal
kernel approximately satisfies the detailed balance condition w.r.t.
$\mu$.

%s1.2.1 #&#
\subsubsection{Ornstein--Uhlenbeck proposals}
A straightforward approach would be to use a proposal density that
satisfies the detailed balance condition
%
%
%e1.20 #&#
\begin{equation}
\gamma^d(x)p(x,y) = \gamma^d(y)p(y,x) \qquad\mbox{for any
}x,y\in\mathbb R^d \label{eq22}
\end{equation}
w.r.t. the standard normal distribution. In this case,
% for $V\equiv0$ all proposed moves are accepted, and in general,
%
%
%e1.21 #&#
\begin{equation}
G(x,y) = V(y)-V(x). \label{eq23}
\end{equation}
The simplest form of proposal distributions satisfying (\ref{eq22}) are
the transition kernels of AR$(1)$ (discrete Ornstein--Uhlenbeck)
processes given by
%
%
%e1.22 #&#
\begin{equation}
p_h^{\mathrm{OU}}(x,dy) = N \biggl( \biggl(1-\frac h2 \biggr)x,
\biggl(h-\frac
{h^2}4 \biggr)I_d \biggr) \label{eq24}
\end{equation}
for some constant $h\in(0,2)$. If $Z$ is a standard normally
distributed $\mathbb R^d$-valued random variable, then the random variables
%
%
%e1.23 #&#
\begin{equation}
\label{eq24a}Y_h^{\mathrm{OU}}(x):= \biggl(1-\frac h2 \biggr)x+
\sqrt{h-\frac{h^2}4}Z,\qquad x\in\mathbb R^d,
\end{equation}
have distributions $p_h^{\mathrm{OU}}(x,dy)$. Note that by (\ref
{eq23}), the acceptance
probabilities
%
%
%e1.24 #&#
\begin{equation}
\label{eq24b}\alpha^{\mathrm{OU}}(x,y) = \exp \bigl(-G^\mathrm{OU}(x,y)^+
\bigr) = \exp \bigl( -\bigl(V(y)-V(x)\bigr)^+ \bigr)
\end{equation}
for Ornstein--Uhlenbeck proposals do not depend on $h$.

%s1.2.2 #&#
\subsubsection{Euler proposals}
In continuous time,
under appropriate regularity and growth conditions on $V$, detailed
balance w.r.t. $\mu$ is satsfied exactly by
the transition functions of the diffusion process solving the
over-damped Langevin stochastic differential equation
%
%
%e1.25 #&#
\begin{equation}
dX_t = -\tfrac12X_t \,dt-\tfrac12\nabla
V(X_t) \,dt+dB_t, \label{eq26}
\end{equation}
because the generator
\[
\mathscr L = \tfrac12\Delta-\tfrac12 x\cdot\nabla-\tfrac12\nabla V\cdot\nabla=
\tfrac12(\Delta-\nabla U\cdot\nabla)
\]
is a self-adjoint operator on an appropriate dense subspace of
$L^2(\mathbb R^d;\mu)$; cf. \cite{Roy}. Although we cannot compute and
sample from the transition functions exactly, we can use approximations
as proposals in a Metropolis--Hastings algorithm. A corresponding MH
algorithm where the proposals are obtained from a discretization scheme
for the SDE (\ref{eq26}) is called a \emph{Metropolis-adjusted Langevin
algorithm} (MALA); cf. \cite{RT,RC}.

In this paper, we focus on the MALA scheme with proposal kernel
%
%
%e1.26 #&#
\begin{equation}
p_h(x,\cdot) = N \biggl( \biggl(1-\frac h2 \biggr)x-\frac h2\nabla
V(x), \biggl(h-\frac{h^2}4 \biggr)\cdot I_d \biggr)
\label{eq31}
\end{equation}
for some constant $h \in(0,2)$; that is, $p_h(x,\cdot)$ is the
distribution of
%
%
%e1.27 #&#
\begin{eqnarray}\label{eq32}
Y_h(x)&=&x-\frac h2\nabla U(x)+\sqrt{h-\frac{h^2}4}Z
\nonumber
\\[-8pt]
\\[-8pt]
\nonumber
&=& \biggl(1-\frac h2 \biggr)x-\frac h2\nabla V(x)+\sqrt{h-\frac{h^2}4}Z ,
\end{eqnarray}
where $Z\sim\gamma^d$ is a standard normal random variable with values
in $\mathbb R^d$.\vadjust{\goodbreak}

Note that if $h-h^2/4$ is replaced by $h$, then (\ref{eq32}) is a
standard Euler discretization step for the
SDE (\ref{eq26}). Replacing $h$ by $h-h^2/4$ ensures that detailed
balance is satisfied exactly for $V\equiv0$.
Alternatively, (\ref{eq32}) can be viewed as a semi-implicit Euler
discretization step for (\ref{eq26}):

%
%re1.2 #&#
\begin{rmk}[(Euler schemes)]\label{ex2}
The \emph{explicit Euler discretization} of the over-damped Langevin
equation (\ref{eq26}) with time step size $h>0$ is given by
%
%
%e1.28 #&#
\begin{equation}\qquad
X_{n+1} = \biggl(1-\frac h2 \biggr)X_n-\frac h2\nabla
V(X_n)+\sqrt h Z_{n+1},\qquad n=0,1,2,\ldots, \label{eq27}
\end{equation}
where $Z_n,n\in\mathbb N$, are i.i.d. random variables with
distribution $\gamma^d$. The process $(X_n)$ defined by (\ref{eq27}) is
a time-homogeneous Markov chain with transition kernel
%
%
%e1.29 #&#
\begin{equation}
p_h^{\mathrm{Euler}}(x,\cdot) = N \biggl( \biggl(1-\frac{h}2
\biggr)x-\frac h 2\nabla V(x),h \cdot I_d \biggr). \label{eq28}
\end{equation}
Even for $V\equiv0$, the measure $\mu$ is not a stationary
distribution for the kernel $p_h^{\mathrm{Euler}}$.
%In fact, if $X$ and $Z$ are independent with distribution $\gamma^d$
%then
%
A \emph{semi-implicit Euler scheme} for (\ref{eq26}) with time-step size
$\epsilon>0$ is given by
%
%
%e1.30 #&#
\begin{equation}
X_{n+1}-X_n = -\frac\epsilon2\cdot\frac{X_{n+1}+X_n}{2}-\frac
\epsilon2\nabla V(X_n)+\sqrt\epsilon Z_{n+1} \label{eq29}
\end{equation}
with $Z_n$ i.i.d. with distribution $\gamma^d$; cf. \cite{HSV}. Note
that the scheme is implicit only in the linear part of the drift but
explicit in $\nabla V$. Solving for $X_{n+1}$ in (\ref{eq29}) and
%yields the equivalent equation
substituting $h=\epsilon/(1+\frac\epsilon4)$ with $h\in(0,2)$
%, i.e., $\epsilon=h/(1-\frac h4)$, we obtain
yields the equivalent equation
%
%
%e1.31 #&#
\begin{equation}
X_{n+1} = \biggl(1-\frac h2 \biggr)X_n-\frac h2\nabla
V(X_n)+\sqrt {h-\frac{h^2}{4}} Z_{n+1}. \label{eq30}
\end{equation}
%
%because $(1-\frac\epsilon4)/(1+\frac\epsilon4)=1-\frac h2$ and $
\end{rmk}

We call the Metropolis--Hastings algorithm with proposal kernel
$p_h(x,\cdot)$ a \emph{semi-implicit MALA scheme with step size $h$}.

%
%pr1.3 #&#
\begin{prop}[(Acceptance probabilities for semi-implicit
MALA)]\label{propE}
Let $V\in C^1(\mathbb R^d)$ and $h\in(0,2)$. Then the acceptance
probabilities for the Metropolis-adjusted Langevin algorithm with
proposal kernels $p_h$ are given by $\alpha_h(x,y)=\exp(-G_h(x,y)^+)$ with
%
%
%e1.32 #&#
\begin{eqnarray}\label{eq33}
\qquad&&G_h(x,y)\nonumber\\
&&\qquad = V(y)-V(x)-\frac{y-x}{2}\cdot\bigl(\nabla
V(y)+\nabla V(x)\bigr)
\\
&&\qquad\quad{}+\frac{h}{8-2h} \bigl[(y+x)\cdot\bigl(\nabla V(y)-\nabla V(x)\bigr)+\bigl|
\nabla V(y)\bigr|^2-\bigl|\nabla V(x)\bigr|^2 \bigr].\nonumber
\end{eqnarray}
For explicit Euler proposals with step size $h >0$, a corresponding
representation holds with
%
%
%e1.33 #&#
\begin{eqnarray}\label{eq34}
G_h^{\mathrm{Euler}}(x,y) &=&V(y)-V(x)-\frac{y-x}{2}\cdot \bigl(
\nabla V(y)+\nabla V(x)\bigr)
\nonumber
\\[-8pt]
\\[-8pt]
\nonumber
&&{}+\frac{h}{8} \bigl[\bigl|y+\nabla V(y)\bigr|^2-\bigl|x+\nabla
V(x)\bigr|^2 \bigr].
\end{eqnarray}
\end{prop}
The proof of the proposition is given in Section~\ref{secRP}
below.

%
%re1.4 #&#
\begin{rmk}
For explicit Euler proposals, the $O(h)$ correction term
%in (\ref{eq33}) vanishes if $V\equiv0$.
%The corresponding correction term for explicit Euler proposals
in (\ref{eq34}) does not vanish for $V\equiv0$.
More significantly, this term goes to infinity as $|y-x|\to\infty$, and
the variance of
$y-x$ w.r.t. the proposal distribution is of order $O(d)$.
\end{rmk}

%s1.3 #&#
\subsection{Bounds for rejection probabilities}
We fix a norm $\| \cdot\|_-$ on $\mathbb R^d$ such that
%
%
%e1.34 #&#
\begin{equation}
\|x\|_- \leq|x| \qquad\mbox{for any }x\in\mathbb R^d. \label{eq37}
\end{equation}
We assume that $V$ is sufficiently smooth w.r.t. $\| \cdot\|_-$ with
derivatives growing at most polynomially:

%
%as1.5 #&#
\begin{assumption}\label{assA1}
The function $V$ is in $C^4(\mathbb R^d)$, and for any $n\in\{
1,2,3,4\}
$, there exist finite constants $C_n\in[0,\infty)$, $p_n\in\{
0,1,2,\ldots\} $ such that
\[
\bigl|\bigl(\partial^n_{\xi_1,\ldots,\xi_n}V\bigr) (x)\bigr| \leq C_n
\max\bigl(1,\|x\| _-\bigr)^{p_n}\| \xi_1\|_-\cdot\cdots\cdot\|
\xi_n\|_-
\]
holds for any $x\in\mathbb R^d$ and $\xi_1,\ldots,\xi_n\in\mathbb R^d$.
\end{assumption}
For discretizations of infinite-dimensional models, $\|\cdot\|_-$ will
typically be a
finite-dimensional approximation of a norm that is almost surely finite
w.r.t. the
limit measure in infinite-dimensions.

%
%ex1.6 #&#
\begin{ex}[(Transition path sampling)]\label{exTPS2}
Consider the situation of Example~\ref{exTPS1}, and assume that $H$ is
in $C^6 (\mathbb R^d)$. Then by (\ref{eqexcirccirc}) and (\ref
{eqexdotdot}), $V_d$ is $C^4$. For $n\le4$
and $x,\xi_1,\ldots,\xi_n\in\mathbb R^d$, the directional derivatives
of $V_d$ are given by
%
%
%e1.35 #&#
\begin{equation}
\label{eqderiv}\quad  \partial^n_{\xi_1\cdots\xi_n}V_d(x) =
2^{-m-1}\sum_{k=0}^{2^m}w_k
D^n\phi (y_{k2^{-m}}) [ \eta_{1,k2^{-m}},\ldots,
\eta_{n,k2^{-m}} ],
\end{equation}
where $y,\eta_1,\ldots,\eta_n$ are the polygonal paths in $E$
corresponding to
$x,\xi_1,\ldots,\xi_n$, respectively, $w_k=1$ for $k=1,\ldots,2^m-1$
and $w_0=w_1=1/2$. Assuming $\| D^4\phi(z)\| =O(|z|^r)$ for some
integer $r\ge0$
as $|z|\to\infty$, we can estimate
\[
\bigl\llvert \partial^n_{\xi_1\cdots\xi_n}V_d(x)\bigr\rrvert
\le C_n \max \bigl(1,\| y\|_{L^q}\bigr)^{p_n} \|
\eta_1\|_{L^q}\cdot\cdots\cdot\|\eta_n
\|_{L^q},
\]
where $q=r+4$, $p_n=r+(4-n)$, $\| y\|_{L^q}=2^{-m}\sum_{k=0}^{2^m}w_k|y_k|^q$
is a discrete $L^q$ norm of the polygonal path $y$ and $C_1,\ldots
,C_4$ are finite
constants that \emph{do not depend on the dimension $d$}. One could now
choose for the minus
norm the norm on $\mathbb R^d$ corresponding to the discrete $L^q$ norm on
polygonal paths. However, it is more convenient to choose a norm coming
from an inner
product. To this end, we consider the norms
\[
\| y\|_\alpha = \biggl( \sum_{n,k,i}2^{-2\alpha n}x_{n,k,i}^2
\biggr)^{1/2},\qquad y=e+\sum x_{n,k,i}e^{n,k,i},
\]
on path space $E$, and the induced norms
\[
\| x\|_\alpha = \biggl( \sum_{n<m}\sum
_{k,i}2^{-2\alpha
n}x_{n,k,i}^2
\biggr)^{1/2}, \qquad x\in\mathbb R^d,
\]
on $\mathbb R^d$ where $d=(2^m-1)\ell$. One can show that for $\alpha<
1/2+1/q$,
the $L^q$ norm can be bounded from above by $\|\cdot\|_\alpha$
independently of the dimension; cf. \cite{G}. On the other hand, if
$\alpha>1/2$, then $\| y\|_\alpha<\infty$
for $\gamma$-almost every path $y$ of the Brownian bridge. This
property will be crucial
when restricting to balls w.r.t. $\|\cdot\|_\alpha$. For $\|\cdot\|
_-=\|\cdot\|_\alpha$ with
$\alpha\in(1/2,1/2+1/q)$, both requirements are satisfied, and
Assumption~\ref{assA1} holds
with constants that do not depend on the dimension.
\end{ex}

The next proposition yields in particular an upper bound for the
average rejection probability w.r.t. both Ornstein--Uhlenbeck and
semi-implicit Euler proposals at a given position $x\in\mathbb R^d$;
cf. \cite{BouRabee} for an analogue result:

%
%pr1.7 #&#
\begin{prop}[(Upper bounds for MH rejection
probabilities)]\label{propI}
Suppose that Assumption~\ref{assA1} is satisfied and let $k\in
\mathbb N$. Then there exist polynomials $\mathcal P_k^{\mathrm
{OU}}\dvtx \mathbb
R\to\mathbb R_+$ and $\mathcal P_k\dvtx \mathbb
R^2\to\mathbb R_+$ of degrees $p_1+1$, $\max(p_3+3,2p_2+2)$,
respectively, such that for any $ x\in\mathbb R^d$ and $h\in(0,2)$,
\begin{eqnarray*}
E\bigl[\bigl(1-\alpha^{\mathrm{OU}}\bigl(x,Y_h^{\mathrm{OU}}(x)
\bigr)\bigr)^k\bigr]^{1/k} &\leq& \mathcal
P_k^{\mathrm{OU}}\bigl(\|x\|_-\bigr)\cdot h^{1/2}\quad\mbox{and}
\\
E\bigl[\bigl(1-\alpha_h\bigl(x,Y_h(x)\bigr)
\bigr)^k\bigr]^{1/k} &\leq&\mathcal P_k\bigl(
\|x\|_-,\bigl\|\nabla U(x)\bigr\|_-\bigr)\cdot h^{3/2}.
\end{eqnarray*}
\end{prop}

The result is a consequence of Proposition~\ref{propE}. The proof is
given in Section~\ref{secRP} below.

%
%re1.8 #&#
\begin{rmk}(1) The polynomials $\mathcal P_k^{\mathrm{OU}}$ and
$\mathcal P_k$ in Proposition~\ref{propI} are explicit; cf. the proof below.
They depend only on the values $C_n, p_n$ in Assumption~\ref{assA1}
for $n=1$,
$n=2,3$, respectively, and on the moments
\begin{eqnarray}
m_n = E\bigl[\|Z\|_-^n\bigr], \qquad n\leq k
\cdot(p_1+1),\qquad  n\leq k\cdot\max(p_3+3,2p_2+2),\nonumber\\
\eqntext{\mbox{respectively},}
\end{eqnarray}
but they do not depend on the dimension $d$. For semi-implicit Euler proposals,
the upper bound in
Proposition~\ref{propI} is stated in explicit form for the case $k=1$
and $p_2=p_3=0$ in (\ref{eqPE}) below.

(2) For \emph{explicit Euler proposals}, corresponding estimates hold
with $m_n$ replaced by $\tilde m_n = E[|Z|^n]$; cf. Remark~\ref
{rmexplicit} below. Note, however, that $\tilde m_n\to\infty$ as
$d\to
\infty$.
%(3) In the case of \emph{Ornstein--Uhlenbeck proposals}, the
%acceptance probabilities of the Metropolis--Hastings algorithm are
%given by (\ref{eq7}) with $G(x,y)=V(y)-V(x)$ independently
%of the step size $h$. From this representation it is not difficult to
%derive a similar bound
%as in Proposition~\ref{propI} but with $h^{3/2}$ replaced by $h^{1/2}$.
\end{rmk}

Our next result is a bound of order $O(h^{1/2})$, $O(h^{3/2})$,
respectively, for
the average
dependence of the acceptance event on the current state w.r.t. Ornstein--Uhlenbeck
and semi-implicit Euler proposals. Let
$\| \cdot\|_+$ denote the dual norm of $\| \cdot\|_-$ on
$\mathbb R^d$, that is,
\[
\|\xi\|_+ = \sup\bigl\{\xi\cdot\eta| \eta\in\mathbb R^d \mbox{ with
} \| \eta\|_-\leq1\bigr\}.
\]
Note that
\[
\|\xi\|_- \leq|\xi| \leq\|\xi\|_+ \qquad \forall\xi\in\mathbb R^d.
\]
For a function $F\in C^1(\mathbb R^d)$,
\begin{eqnarray*}
\bigl|F(y)-F(x)\bigr|&=&\biggl\llvert \int_0^1(y-x)
\cdot\nabla F\bigl((1-t)x+ty\bigr) \,dt\biggr\rrvert
\\
&\leq&\|y-x\|_-\cdot\sup_{z\in[x,y]}\bigl\|\nabla F(z)\bigr\|_+,
\end{eqnarray*}
that is, the plus norm of $\nabla F$ determines the Lipschitz constant
w.r.t. the minus norm.

%
%pr1.9 #&#
\begin{prop}[(Dependence of rejection on the current
state)]\label{propM}
Suppose that Assumption~\ref{assA1} is satisfied, and let
$k\in\mathbb N$. Then there exist polynomials $\mathcal Q_k^{\mathrm
{OU}}\dvtx \mathbb
R\to\mathbb R_+$ and $\mathcal Q_k\dvtx \mathbb
R^2\to\mathbb R_+$ of degrees $p_2+1$, $\max(p_4+3,p_3+p_2+2,3p_2+1)$,
respectively, such that for any $ x,\tilde x\in\mathbb R^d$ and $h\in(0,2)$,
%
%
%e1.36 #&#
%e1.37 #&#
\begin{eqnarray}\qquad
&&E\bigl[\bigl\|\nabla_xG^\mathrm{OU}\bigl(x,Y_h^{\mathrm{OU}}(x)
\bigr)\bigr\| _+^k\bigr]^{1/k}\leq \mathcal
Q_k^{\mathrm{OU}} \bigl(\|x\|_- \bigr)\cdot h^{1/2}, \label{eqMa}
\\
&&E\bigl[\bigl\|\nabla_xG_h\bigl(x,Y_h(x)\bigr)
\bigr\|_+^k\bigr]^{1/k}\leq\mathcal Q_k \bigl(\|x
\| _-,\bigl\|\nabla U(x)\bigr\|_- \bigr)\cdot h^{3/2},\label{eqMb}\\
%%
%%
%%e1.38 #&#
%%e1.39 #&#
\label{eqMc}&&E\bigl[\bigl|\alpha^{\mathrm{OU}}\bigl(x,Y_h^{\mathrm{OU}}(x)
\bigr)-\alpha ^{\mathrm
{OU}}\bigl(\tilde x,Y_h^{\mathrm{OU}}(\tilde
x)\bigr)\bigr|^k\bigr]^{1/k}
\nonumber
\\[-8pt]
\\[-8pt]
\nonumber
&&\qquad \leq\mathcal Q_k^{\mathrm{OU}} \bigl(\max\bigl(\|x\|_-,\|\tilde x\|
_-\bigr) \bigr)\cdot \|x-\tilde x\|_-\cdot h^{1/2}\quad \mbox{and}
\\
\label{eqMd}
&&E\bigl[\bigl|\alpha_h\bigl(x,Y_h(x)\bigr)-
\alpha_h\bigl(\tilde x,Y_h(\tilde x)
\bigr)\bigr|^k\bigr]^{1/k}
\nonumber
\\[-8pt]
\\[-8pt]
\nonumber
&&\qquad \leq\mathcal Q_k \Bigl(\max\bigl(\|x\|_-,\|\tilde x\|_-\bigr),\sup
_{z\in
[x,\tilde
x]}\bigl\|\nabla U(z)\bigr\|_- \Bigr)\cdot\|x-\tilde x\|_-\cdot
h^{3/2},
\end{eqnarray}
where $[x,\tilde x]$ denotes the line segment between $x$ and
$\tilde x$.
\end{prop}
The proof of the proposition is given in Section~\ref{secDRP}
below.\vadjust{\goodbreak}

%
%re1.10 #&#
\begin{rmk}
Again, the polynomials $\mathcal Q_k^{\mathrm{OU}}$ and $\mathcal Q_k$
are explicit.
They depend only on the values $C_n, p_n$ in Assumption~\ref{assA1}
for $n=1,2$,
$n=2,3,4$, respectively, and on the moments
$m_n=E[\|Z\|_-^n]$ for $n\leq k\cdot(p_2+1)$, $n\leq k\cdot\max
(p_4+3,p_3+p_2+2,2p_2+1)$, respectively,
but they do not depend on the dimension $d$. For semi-implicit Euler proposals,
the upper bound in
Proposition~\ref{propM} is made explicit for the case $k=1$ and
$p_2=p_3=p_4=0$ in (\ref{eqQE}) below.
\end{rmk}

For Ornstein--Uhlenbeck proposals, it will be useful to state the
bounds in
Propositions \ref{propI} and \ref{propM} more explicitly for the case
$p_2=0$, that is,
when the
second derivatives of $V$ are uniformly bounded w.r.t. the minus norm:

%
%pr1.11 #&#
\begin{prop}\label{propMOU}
Suppose that Assumption~\ref{assA1} is satisfied for $n=1, 2$ with
$p_2=0$. Then for any $ x,\tilde x\in\mathbb R^d$ and $h\in(0,2)$,
\begin{eqnarray*}
&&\mathbb{E}\bigl[1-\alpha^{\mathrm{OU}}\bigl(x,Y_h^{\mathrm{OU}}(x)
\bigr)\bigr] \\
&&\qquad\le m_1 \bigl(C_1+C_2 \| x\|_-\bigr)\cdot
h^{1/2}
\\
& &\qquad\quad{}+ \tfrac12 \bigl(2m_2C_2+C_1\| x
\|_-+C_2\| x\|_-^2\bigr)\cdot h + \tfrac12m_1C_2
\| x\|_-\cdot h^{3/2}, %E[\|\nabla_xG(x,Y_h^{\mathrm{OU}}(x))\|_+^k]^{1/k}\ \leq\ C_2
%m_k^{1/k}\cdot h^{1/2} + (C_2\| x\|_-+C_1/2)\cdot h
\end{eqnarray*}
and
\begin{eqnarray*}
&&E\bigl[\bigl|\alpha^{\mathrm{OU}}\bigl(x,Y_h^{\mathrm{OU}}(x)
\bigr)-\alpha ^{\mathrm
{OU}}\bigl(\tilde x,Y_h^{\mathrm{OU}}(\tilde
x)\bigr)\bigr|^2\bigr]^{1/2}
\\
&&\qquad\leq \bigl( m_2^{1/2} C_2\cdot
h^{1/2} + \tfrac12 \bigl(C_1+2C_2\max \bigl(\| x\|_-,
\|\tilde x\|_- \bigr) \bigr)\cdot h \bigr)\cdot\| x-\tilde x\|_-.
\end{eqnarray*}
\end{prop}

The proof is given in Sections \ref{secRP} and \ref{secDRP} below. Again,
corresponding bounds also hold for $L^k$ norms for $k\neq1,2$.
%We now endow $\mathbb R^d$ with the bounded metrics
%The centered ball of radius $R$ w.r.t. $\| \cdot\|_-$ is denoted by
%d_R(x,0)\neq R\}.\]
%We assume that $V$ is in $C^4(\mathbb R^d)$ and Assumption~\ref{assA1}
%holds. Let $P$ and $Q$ denote the %polynomials from Proposition~\ref{propI} and Proposition~\ref{propM} respectively.

%s1.4 #&#
\subsection{Wasserstein contractivity}
The bounds in Propositions \ref{propI}, \ref{propM} and \ref{propMOU}
can be applied to study contractivity properties of
Metropolis--Hastings transition kernels.
Recall that the \emph{Kantorovich--Rubinstein} or \emph
{$L^1$-Wasserstein distance} of two probability
measures $\mu$ and $\nu$ on the Borel $\sigma$-algebra $\mathcal
B(\mathbb R^d)$ w.r.t. a given metric~$d$ on $\mathbb R^d$ is defined by
\[
\mathcal W(\mu,\nu) = \inf_{\eta\in\Pi(\mu,\nu)}\int d(x,\tilde x) \eta(dx\, d
\tilde x),
\]
where $\Pi(\mu,\nu)$ consists of all couplings $\eta$ of $\mu$
and $\nu$, that is, all probability measures $\eta$ on $\mathbb
R^d\times\mathbb R^d$ with marginals $\mu$ and $\nu$;
cf., for example, \cite{Villani}. Recall that a coupling of $\mu$ and
$\nu$ can be realized by random variables $W$ and $\tilde W$ defined on
a joint probability space such that $W\sim\mu$
and $\tilde W\sim\nu$.

In order to derive upper bounds for the
distances $\mathcal W(\mu q_h,\nu q_h)$, and, more
generally, $\mathcal W(\mu q_h^n,\nu q_h^n),n\in\mathbb N$, we
define a coupling of the MALA transition probabilities
$q_h(x, \cdot),x\in\mathbb R^d$, by setting
\[
W_h(x):= %
\cases{Y_h(x),&\quad$\mbox{if } {
\mathcal U}\leq\alpha_h\bigl(x,Y_h(x)\bigr),$
\vspace*{2pt}
\cr
x,&$\quad\mbox{if } { \mathcal U}>\alpha_h
\bigl(x,Y_h(x)\bigr).$} %
\]
Here $Y_h(x)$, $x\in\mathbb R^d$, is the basic coupling of the
proposal distributions $p_h(x, \cdot)$ defined by
(\ref{eq32}) with $Z\sim\gamma^d$, and the random variable
${\mathcal
U}$ is uniformly distributed in $(0,1)$ and independent of $Z$.

Correspondingly, we define a coupling
of the Metropolis--Hastings transition kernels $q_h^\mathrm{OU}$ based
on Ornstein--Uhlenbeck proposals by setting
\[
W_h^\mathrm{OU}(x):= %
\cases{Y_h^\mathrm{OU}(x),&\quad $
\mbox{if } { \mathcal U}\leq\alpha ^\mathrm{OU}\bigl(x,Y_h^\mathrm{OU}(x)
\bigr),$\vspace*{2pt}
\cr
x,&\quad $\mbox{if } { \mathcal U}>\alpha ^\mathrm{OU}
\bigl(x,Y_h^\mathrm{OU}(x)\bigr).$} %
\]

Let
\[
B_R^-:= \bigl\{x\in\mathbb R^d\dvtx \|x\|_-<R\bigr\}
\]
denote the centered ball of radius $R$ w.r.t. $\| \cdot\|_-$. As a
consequence of
Proposition~\ref{propMOU} above, we obtain the following upper bound
for the
Kantorovich--Rubinstein--Wasserstein distance of $q_h^\mathrm{OU}(x,
\cdot)$ and $q_h^\mathrm{OU}(\tilde x, \cdot)$ w.r.t. the
metric $d(x,\tilde x)=\| x-\tilde x\|_-$:

%
%th1.12 #&#
\begin{thmm}[(Contractivity of MH
transitions based on OU proposals)]\label{thmQOU}
Suppose that Assumption~\ref{assA1} is satisfied for $n=1, 2$ with
$p_2=0$. Then for any $h\in(0,2)$, $R\in(0,\infty)$, and $x,\tilde
x\in B_ R^-$,
\[
\mathbb{E}\bigl[\bigl\| W_h^\mathrm{OU}(x)-W_h^\mathrm{OU}(
\tilde x)\bigr\|_- \bigr] \leq c_h^\mathrm{OU}(R)\cdot\|x-\tilde x
\|_-,
\]
where
\[
c_h^\mathrm{OU}(R) = 1 - \tfrac12 h + m_2C_2
h + A(1+R) \bigl(1+h^{1/2}R\bigr) h^{3/2}
\]
with an explicit constant $A$ that only depends on the values
$m_1,m_2,C_1$ and~$C_2$.
\end{thmm}
The proof is given in Section~\ref{secFinal}
below.

Theorem~\ref{thmQOU} shows that Wasserstein contractivity holds on the
ball $B_R^-$
provided $2m_2C_2<1$ and $h$ is chosen sufficiently small depending on $R$
[with $h^{1/2}=O (R^{-1})$]. In this case, the contraction constant
$c_h^\mathrm{OU}(R)$ depends on the dimension only through the values
of the constants $C_1, C_2, m_1$ and $m_2$. On the other hand, the
following one-dimensional example shows that for $m_2C_2>1$, the
acceptance-rejection step may destroy the contraction properties of the
OU proposals:

%
%ex1.13 #&#
\begin{ex}Suppose that $d=1$ and $\|\cdot\|_-=|\cdot| $. If
$V(x)=bx^2/2$ with
a constant $b\in(-1/2,1/2)$, then by Theorem~\ref{thmQOU},
Wasserstein contractivity holds
for the Metropolis--Hastings chain with Ornstein--Uhlenbeck proposals
on the interval $(-R,R)$ provided $h$ is chosen sufficiently small. On
the other hand, if $V(x)=bx^2/2$ for $|x|\le1$ with a constant $b<-1$,
then the logarithmic density
\[
U(x) = V(x)+x^2/2 = (b+1)\cdot x^2/2
\]
is strictly concave for $|x|\le1 $, and it can
be easily seen that Wasserstein contractivity on $(-1,1)$ does not hold
for the MH chain with
OU proposals if $h$ is sufficiently small.
\end{ex}

A disadvantage of the result for Ornstein--Uhlenbeck proposals stated
above is that
not only a lower bound on the second derivative of $V$ is required
(this would be
a fairly natural condition as the example indicates), but also an upper
bound of the same
size. For semi-implicit Euler proposals, we can derive a better result
that requires only
a strictly positive lower bound on the second
derivative of $U(x)=V(x)+|x|^2/2$ and Assumption~\ref{assA1} with
arbitrary constants
to be satisfied. For this purpose we assume that
\[
\| \cdot\|_- = \langle\cdot, \cdot\rangle^{1/2}
\]
for
an inner product $\langle\cdot, \cdot\rangle$ on $\mathbb
R^d$, and we make the
following assumption on~$U$:

%
%as1.14 #&#
\begin{assumption}\label{assA2}
There exists a strictly positive constant $K\in(0,1]$ such that
%
%
%e1.40 #&#
\begin{equation}
\label{eq60} \bigl\langle\eta,\nabla^2U(x)\cdot\eta\bigr\rangle\geq K
\langle\eta,\eta\rangle\qquad\mbox{for any }x, \eta\in\mathbb R^d.
\end{equation}
\end{assumption}

Of course, Assumption~\ref{assA2} is still restrictive, and it will
often be satisfied only in a suitable
ball around a local minimum of $U$. Most of the results below are
stated on a given
ball $B_R^-$ w.r.t.  the minus norm. In this case it is enough to
assume that \ref{assA2}
holds on that ball. If $\| \cdot\|_-$ coincides with the Euclidean
norm $| \cdot|$,
then the assumption is equivalent to convexity of $U(x)-K|x|^2$.
Moreover, since $\nabla^2U(x)=I_d+\nabla^2V(x)$, a sufficient
condition for
(\ref{eq60}) to hold is
%
%
%e1.41 #&#
\begin{equation}
\label{eqA2suff} \bigl\llVert \nabla^2V(x)\cdot\eta\bigr\rrVert _-
\leq(1-K) \llVert \eta\rrVert _-\qquad \mbox{for any }x, \eta\in\mathbb
R^d.
\end{equation}

As a consequence of
Propositions \ref{propI} and \ref{propM} above, we obtain the following
upper bound for the
Kantorovich--Rubinstein--Wasserstein distance of $q_h(x, \cdot)$ and
$q_h(\tilde x, \cdot)$ w.r.t. the
metric $d(x,\tilde x )=\| x-\tilde x\|_-$:

%
%th1.15 #&#
\begin{thmm}[(Contractivity of semi-implicit MALA
transitions)]\label{thmQ}
Suppose that Assumptions \ref{assA1} and \ref{assA2} are satisfied.
Then for any $h\in(0,2)$, $R\in(0,\infty)$ and $x,\tilde x\in B_ R^-$,
\[
\mathbb{E} \bigl[\bigl\| W_h(x)-W_h(\tilde x)\bigr\|_- \bigr]
\leq c_h(R)\cdot \| x-\tilde x\|_-,
\]
where
\[
c_h(R) = 1 - \tfrac12 K h + \bigl(\tfrac18 M(R)^2+
\gamma(R) \bigr) h^{2} + \bigl( K\beta(R)+\tfrac12\delta(R) \bigr)
h^{5/2}
\]
with
\begin{eqnarray*}
M(R)&=&\sup\bigl\{\bigl\|\nabla^2U(z)\cdot\eta\bigr\|_-\dvtx \eta\in
B_1^-, z\in B_R^-\bigr\} ,
\\
\beta(R)&=&\sup\bigl\{\mathcal P_1\bigl(\|z\|_-,\bigl\|\nabla U(z)\bigr\|_-
\bigr)\dvtx z\in B_R^-\bigr\} ,
\\
\gamma(R)&=& m_2^{1/2}\cdot\sup\bigl\{\mathcal
Q_2\bigl(\|z\|_-,\bigl\|\nabla U(z)\bigr\| _-\bigr)\dvtx z\in B_R^-
\bigr\},
\\
\delta(R)&=&\sup\bigl\{\mathcal Q_2\bigl(\|z\|_-,\bigl\|\nabla U(z)\bigr\|_-
\bigr) \bigl\|\nabla U(z)\bigr\|_-\dvtx z\in B_R^-\bigr\}.
\end{eqnarray*}
\end{thmm}

The proof is given in Section~\ref{secFinal}
below.

%
%re1.16 #&#
\begin{rmk}\label{rmkQ}
Theorem~\ref{thmQ} shows in particular that under Assumptions \ref
{assA1} and \ref{assA2}, there exist constants $C,q\in(0,\infty)$
such that the contraction
\[
\mathbb{E}\bigl[\bigl\| W_h(x)-W_h(\tilde x)\bigr\|_- \bigr]
\leq \biggl( 1-\frac K4 h \biggr) \| x-\tilde x\|
\]
holds for $x,\tilde x\in B_R^-$ whenever $h^{-1}\ge C\cdot(1+R^q)$.
\end{rmk}

%
%ex1.17 #&#
\begin{ex}[(Transition path sampling)]
In the situation of Examples \ref{exTPS1} and \ref{exTPS2} above, condition
(\ref{eqA2suff}) and (hence) assumption \ref{assA2} are satisfied on
a ball
$B_R^-$ with $K$ independent of $d$ provided $\| D^2\phi(x)\|\le1-K$
for any $x\in B_R^-$; cf.~(\ref{eqderiv}). More generally, by
modifying the metric in
a suitable way if necessary, one may expect Assumption~\ref{assA2} to hold
uniformly in the dimension in neighborhoods of local minima of $U$.
\end{ex}

%The upper bound in Theorem~\ref{thmQ} can be slightly improved if
%$\| \cdot\|_- = \langle\cdot, \cdot\rangle^{1/2}$ for
%an inner product $\langle\cdot, \cdot\rangle$ on $\mathbb
%R^d$, and there exists a constant $K (R)\in\mathbb R$ s.t.
%In this case,
%&&+\left(\beta(R)h^{3/2}+\gamma(R)h^{2}+\delta
%(R)h^{5/2}\right)\cdot\|x-\tilde x\|_-,
%cf. Section~\ref{secFinal} below.

%s1.5 #&#
\subsection{Conclusions}

For $R\in(0,\infty)$, we denote by $\mathcal W_R$ the
Kantorovich--Rubinstein--Wasserstein distance based on the distance function
%
%
%e1.42 #&#
\begin{equation}
\label{eqdr} d_R(x,\tilde x ) := \min\bigl(\| x-\tilde x\|_-, 2R\bigr).
\end{equation}
Note that $d_R$ is a bounded metric that coincides with the distance
function induced by the minus norm on $B_R^-$.
The bounds resulting from
Theorems \ref{thmQ} and \ref{thmQOU} can be iterated to obtain
estimates for the
KRW distance $\mathcal W_R$ between the
distributions of the corresponding Metropolis--Hastings chains after
$n$ steps w.r.t. two
different initial distributions.

%
%co1.18 #&#
\begin{coroll}\label{corWbound}
Suppose that Assumptions \ref{assA1} and \ref{assA2} are satisfied, and
let $h\in(0,2)$ and $R\in(0,\infty)$. Then for any $n\in\mathbb N$,
and for any probability measures $\mu,\nu$ on $\mathcal B(\mathbb R^d)$,
\begin{eqnarray*}\label{eqWbound}
\mathcal W_R\bigl(\mu q_h^n,\nu
q_h^n\bigr)&\leq&c_h(R)^n
\mathcal W_R(\mu,\nu )
\\
&&{}+ {2R}\cdot \bigl(\P_\mu\bigl[\exists k<n\dvtx X_k\notin
B_R^-\bigr] +\P_\nu\bigl[\exists k<n\dvtx X_k
\notin B_R^-\bigr]\bigr).
\end{eqnarray*}
Here $c_h(R)$ is the constant in Theorem~\ref{thmQ}, and
$(X_n,\P_\mu)$ and $(X_n,\P_\nu)$ are Markov chains with
transition kernel $q_h$ and initial distributions $\mu$, $\nu$,
respectively. A corresponding result with $c_h$ replaced by
$c_h^\mathrm
{OU}$ holds for the Metropolis--Hastings chain with Ornstein--Uhlenbeck
proposals.
\end{coroll}

Since the joint law of $W_h(x)$ and $W_h(\tilde x)$ is a coupling of
$q_h(x,\cdot)$ and
$q_h(\tilde x,\cdot)$ for any $x,\tilde x\in\mathbb R^d$, Corollary~\ref{corWbound} is a
direct consequence of Theorems \ref{thmQ}, \ref{thmQOU}, respectively, and
Theorem~\ref{thmC} below. The corollary can be used to quantify the
Wasserstein distance between
the distribution of the Metropolis--Hastings chain after $n$ steps
w.r.t. two different initial
distributions. For this purpose, one can estimate the exit
probabilities from the
ball $B_R^-$ via an argument based on a Lyapunov function. For
semi-implicit Euler proposals we
eventually obtain the following main result:

%
%th1.19 #&#
\begin{thmm}[(Quantitative convergence bound for semi-implicit
MALA)]\label{thmMAIN}
Suppose that Assumptions \ref{assA1} and \ref{assA2} are satisfied.
Then there exist
constants $C,D,q\in(0,\infty)$ such that the estimate
\[
\mathcal W_{2R}\bigl(\nu q_h^n,\pi
q_h^n\bigr) \leq \biggl( 1-\frac K4 h
\biggr)^n \mathcal W_{2R}(\nu,\pi) + {D R} \exp \biggl(-
\frac{KR^2}{8} \biggr) nh
\]
holds for any $n\in\mathbb N$, $h, R\in(0,\infty)$ such that
$h^{-1}\ge C\cdot(1+R)^q$, and for any probability measures $\nu,\pi
$ on $\mathbb R^d$ with support in $B_{R}^-$. The constants~$C$, $D$
and $q$ can be made explicit. They depend
only on the values of the constants in Assumptions~\ref{assA1} and
\ref
{assA2} and on the
moments $m_k$, $k\in\mathbb N$, w.r.t. the minus norm, but they do not
depend explicitly on the dimension.
\end{thmm}

The proof of Theorem~\ref{thmMAIN} is given in Section~\ref{secFinal}
below.

Let $\mu_R(A)=\mu(A|B_R^-)$ denote the conditional probability measure
given $B_R^-$.
Recalling that $\mu$ is a stationary distribution for the kernel $q_h$,
we can apply Theorem~\ref{thmMAIN} to derive a bound for the Wasserstein distance of the
discretization of the
MALA chain and $\mu_R$ after $n$ steps:

%
%th1.20 #&#
\begin{thmm}\label{thmFINAL}
Suppose that Assumptions \ref{assA1} and \ref{assA2} are satisfied.
Then there exist
constants $C,\bar D,q\in(0,\infty)$ that do not depend explicitly on
the dimension such that the estimate
\[
\mathcal W_{2R}\bigl(\nu q_h^n,
\mu_R \bigr) \leq 58 R \biggl( 1-\frac K4 h \biggr)^n + {
\bar D R} \exp \biggl(-\frac{KR^2}{33} \biggr) nh
\]
holds for any $n\in\mathbb N$, $h, R\in(0,\infty)$ such that
$h^{-1}\ge C\cdot(1+R)^q$, and for any probability measure $\nu$ on
$\mathbb R^d$ with support in $B_{R}^-$.
\end{thmm}

The proof is given in Section~\ref{secFinal}.

Given an error bound $\varepsilon\in(0,\infty)$ for the
Kantorovich--Rubinstein--Wasser\-stein
distance, we can now determine how many steps of the MALA chain are
required such
that
%
%
%e1.43 #&#
\begin{equation}
\label{eq7o} \mathcal W_{2R}\bigl(\nu q_h^n,
\mu_R \bigr) < \varepsilon\qquad \mbox{for any }\nu \mbox{ with support in
}B_R^-.
\end{equation}
Assuming
%
%
%e1.44 #&#
\begin{equation}
\label{eq7star} nh \ge \frac4K \log \biggl( \frac{116R}{\varepsilon} \biggr),
\end{equation}
we have $58 R (1-Kh/4)^n\le\varepsilon/2$. Hence (\ref{eq7o}) holds
provided the
assumptions in Theorem~\ref{thmFINAL} are satisfied, and
%
%
%e1.45 #&#
\begin{equation}
\label{eq7starstar} \bar DR\exp\bigl(-KR^2/33\bigr) nh <
\varepsilon/2.
\end{equation}
For a minimal choice of $n$, all conditions are satisfied if $R$ is of
order $({\log\varepsilon^{-1}})^{1/2}$ up to a $\log\log$ correction,
and the inverse step size
$h^{-1}$ is of order $(\log\varepsilon^{-1})^{q/2}$ up to a $\log
\log$
correction.
Hence if Assumption~\ref{assA2} holds on $\mathbb R^d$, then a number
$n$ of
steps that is polynomial in $\log\varepsilon^{-1}$ is sufficient to
bound the error by
$\varepsilon$ independently of the dimension.

On the other hand, if Assumption~\ref{assA2} is satisfied only on a
ball $B_R^-$ of
given radius $R$, then a given error bound $\varepsilon$ is definitely
achieved only
provided (\ref{eq7starstar}) holds with the minimal choice for $nh$
satisfying (\ref{eq7star}),
that is, if
%
%
%e1.46 #&#
\begin{equation}
\label{eq7starstarstar} 8\bar DK^{-1}\log\bigl(116 R
\varepsilon^{-1}\bigr) R \exp\bigl(-KR^2/33\bigr) <
\varepsilon.
\end{equation}
If $\varepsilon$ is chosen smaller, then the chain may leave the ball
$B_R^-$ before
sufficient mixing on $B_R^-$ has taken place.

%s2 #&#
\section{Wasserstein contractivity of Metropolis--Hastings
kernels}
\label{secCWC}
\label{secWCMH}

In this section, we first consider an arbitrary stochastic kernel
$q\dvtx S\times\mathcal
B(S)\to[0,1]$ on a metric space $(S,d)$.
Further below, we will choose $S=\mathbb R^d$ and $d(x,y)=\|x-y\|
_-\wedge R$ for some constant $R\in(0,\infty]$, and we will assume
that $q$ is the transition kernel of
a Metropolis--Hastings chain.

The \emph{Kantorovich--Rubinstein} or \emph{$L^1$-Wasserstein distance}
of two probability measures $\mu$ and $\nu$ on the Borel-$\sigma
$-algebra $\mathcal B(S)$ w.r.t. the metric $d$ is defined by
\[
\mathcal W_d(\mu,\nu) = \inf_{\eta}\int d(x,\tilde
x) \eta (dx\,d\tilde x),
\]
where the infimum is over all couplings $\eta$ of $\mu$ and $\nu$,
that is, over all probability measures $\eta$ on $S\times S$ with
marginals $\mu$ and $\nu$; cf., for example, \cite{Villani}. In order
to derive
upper bounds for the Kantorovich distances $\mathcal W_d(\mu q,\nu
q)$, and more generally, $\mathcal W_d(\mu q^n,\nu q^n),n\in\mathbb
N$, we construct couplings between the measures $q(x,\cdot)$ for
$x\in S$, and we derive bounds for the distances $\mathcal
W_d(q(x,\cdot),q(\tilde x,\cdot)),x,\tilde x\in S$.

%
%de2.1 #&#
\begin{defi}
A \emph{Markovian coupling} of the probability measures $q(x,\cdot)$,
$x\in S$, is a stochastic kernel
$c$ on the product space $(S\times S, \mathcal B(S\times S))$ such that
for any $x,\tilde x\in S$, the
distribution of the first and second component under $c((x,\tilde x),dy
\,d\tilde y)$ is $q(x,dy)$ and
$q(\tilde x,d\tilde y)$, respectively.
\end{defi}

%
%ex2.2 #&#
\begin{ex}
(1) Suppose that $(\Omega,\mathcal A,\P)$ is a probability space, and
let $(x,\tilde x,\omega)\mapsto
Y(x,\tilde x)(\omega)$, $(x,\tilde x,\omega)\mapsto\tilde Y(x,\tilde
x)(\omega)$ be product
measurable functions from $S\times S\times\Omega$ to $S$ such that
$Y(x,\tilde x)\sim q(x,\cdot)$ and
$\tilde Y(x,\tilde x)\sim q(\tilde x,\cdot)$ w.r.t. $\P$ for any
$x,\tilde x\in S$. Then the joint
distributions
\[
c\bigl((x,\tilde x),\cdot\bigr) = \P\circ\bigl(Y(x,\tilde x),\tilde Y(x,\tilde
x)\bigr)^{-1},\qquad x,\tilde x\in S,
\]
define a Markovian coupling of the measures $q(x,\cdot),x\in
S$.

(2) In particular, if $(x,\omega)\mapsto Y(x)(\omega)$ is a product
measurable function from $S\times\Omega$ to $S$ such that $Y(x)\sim
q(x,\cdot)$ for any $x\in S$, then
\[
c\bigl((x,\tilde x),\cdot\bigr) = \P\circ\bigl(Y(x),Y(\tilde x)
\bigr)^{-1}
\]
is a Markovian coupling of the measures $q(x,\cdot), x\in S$.
\end{ex}

%Further below, we will only consider couplings of the second type.
Suppose that $(X_n,\tilde X_n)$ on $(\Omega,\mathcal A,\P)$ is a Markov
chain with values in $S\times S$ and
transition kernel $c$, where $c$ is a Markovian coupling w.r.t. the
kernel $q$. Then the components $(X_n)$ and $(\tilde X_n)$ are Markov
chains with transition kernel $q$ and initial distributions given by
the marginals of the initial distribution
of $(X_n,\tilde X_n)$, that is, $(X_n,\tilde X_n)$ is a coupling of
these Markov chains.
We will apply the following general theorem to quantify the deviation
from equilibrium after $n$ steps of the
Markov chain with transition kernel $q$:

%
%th2.3 #&#
\begin{thmm}\label{thmC}
Let $\gamma\in(0,1)$, and let $c((x,\tilde x),dy \,d\tilde y)$ be a
Markovian coupling of the probability measures $q(x,\cdot), x\in S$.
Suppose that $\mathcal O$ is an open subset of $S$, and
assume that the metric $d$ is bounded. Let
\[
\Delta := \operatorname{diam} S = \sup\bigl\{ d(x,\tilde x)\dvtx x,\tilde x \in S\bigr\}.
\]
If the contractivity condition
%
%
%e2.1 #&#
\begin{equation}
\int d(y,\tilde y) c\bigl((x,\tilde x),dy \,d\tilde y\bigr) \leq\gamma\cdot d(x,
\tilde x) \label{eq10}
\end{equation}
holds for any $x,\tilde x\in\mathcal O$, then
%
%
%e2.2 #&#
\begin{eqnarray}\label{eq11}
\qquad&&\mathcal W_d\bigl(\mu q^n,\nu q^n\bigr)
\nonumber
\\[-8pt]
\\[-8pt]
\nonumber
&&\qquad\leq
\gamma^n\mathcal W_d(\mu,\nu )+{\Delta}\cdot\bigl(
\P_\mu[\exists k<n\dvtx X_k\notin\mathcal O ] +
\P_\nu[\exists k<n\dvtx X_k\notin\mathcal O ]\bigr)
\end{eqnarray}
for any $n\in\mathbb N$ and for any probability measures $\mu,\nu$ on
$\mathcal B(S)$. Here
$(X_n,\P_\mu)$ and $(X_n,\P_\nu)$ are Markov chains with transition
kernel $q$ and initial distributions $\mu$,
$\nu$, respectively.
\end{thmm}

\begin{pf}%{Proof of Theorem~\ref{thmC}}
Suppose that $\mu$ and $\nu$ are probability measures on $\mathcal
B(S)$ and $\eta(dx \,d\tilde x)$ is a coupling of $\mu$ and $\nu$. We
consider the coupling chain $(X_n,\tilde X_n)$ on $(\Omega,\mathcal
A,\P
)$ with initial
distribution $\eta$ and transition kernel $c$. Since $(X_n)$ and
$(\tilde X_n)$ are Markov chains with transition
kernel $q$ and initial distributions $\mu$ and $\nu$, we have $\P
\circ
X_n^{-1}=\mu q^n$ and $\P\circ\tilde X_n^{-1}=\nu q^n$ for any $n\in
\mathbb N$. Moreover, by~(\ref{eq10}),
\begin{eqnarray*}
&&\mathbb{E} \bigl[ d(X_n,\tilde X_n);
(X_k,\tilde X_k)\in \mathcal O\times\mathcal O\ \forall
k<n \bigr]
\\
&&\qquad=\mathbb{E} \biggl[ \int d(x_n,\tilde x_n) c
\bigl((X_{n-1},\tilde X_{n-1}),dx_n \,d\tilde
x_n\bigr); (X_k,\tilde X_k)\in\mathcal O
\times \mathcal O\ \forall k<n \biggr]
\\
&&\qquad\le\gamma\mathbb{E} \bigl[ d(X_{n-1},\tilde X_{n-1});
(X_k,\tilde X_k)\in\mathcal O\times\mathcal O\ \forall
k<n-1 \bigr].
\end{eqnarray*}
Therefore, by induction,
\begin{eqnarray*}
\mathcal W_d\bigl(\mu q^n,\nu q^n
\bigr) &\le& \mathbb{E} \bigl[ d(X_n,\tilde X_n) \bigr]
\\
&=&\mathbb{E} \bigl[ d(X_n,\tilde X_n);
(X_k,\tilde X_k)\in\mathcal O\times\mathcal O\ \forall
k<n \bigr]
\\
&&{}+ \mathbb{E} \bigl[ d(X_n,\tilde X_n); \exists k<n\dvtx (X_k,\tilde X_k)\notin\mathcal O\times \mathcal O
\bigr]
\\
&\le&\gamma^n d(x,\tilde x) + \Delta\cdot\P \bigl[ \exists k<n\dvtx (X_k,\tilde X_k)\notin\mathcal O\times\mathcal O
\bigr],
\end{eqnarray*}
which implies (\ref{eq11}).
\end{pf}

%(Further below, when applying the theorem we will choose $S=\mathbb
%R^d$ and
%In this case, the bound (\ref{eq12}) becomes
%with\[C_n(r,\mu):= C_n(B(0,r),\mu) = \sup_{0\leq k<n}(\mu
%q^k)(d(x,0)>\Delta/2).)\]

%
%re2.4 #&#
\begin{rmk}
Theorem~\ref{thmC} may also be useful for studying local equilibration
of a Markov chain within a metastable state. In fact, if $\mathcal O$
is a region of the state space where the process stays with high
probability for a long time, and if a contractivity condition holds on
$\mathcal O$, then the result can be used to bound the
Kantorovich--Rubinstein--Wasserstein distance between the distribution
after a finite number of steps and the stationary distribution
conditioned to $\mathcal O$.
\end{rmk}

From now on, we assume that we are given a Markovian coupling of the
proposal distributions
$p(x,\cdot),x\in\mathbb R^d$, of a Metropolis--Hastings algorithm which
is realized by product
measurable functions $(x,\tilde x,\omega)\mapsto Y(x,\tilde x)(\omega
)$, $\tilde Y(x,\tilde
x)(\omega)$ on a probability space $(\Omega,\mathcal A,\P)$ such that
\[
Y(x,\tilde x) \sim p(x,\cdot) \quad\mbox{and}\quad \tilde Y(x,\tilde x) \sim p(\tilde x,
\cdot)\qquad \mbox{for any }x,\tilde x\in\mathbb R^d.
\]
Let $\alpha(x,y)$ and $q(x,dy)$ again denote the acceptance
probabilities and the transition kernel of the Metropolis--Hastings
chain with stationary distribution $\mu$; cf. (\ref{eq2}), (\ref{eq3})
and (\ref{eq4}). Moreover, suppose that $\mathcal U$ is a uniformly
distributed random variable with values in $(0,1)$ that is independent
of $\{Y(x,\tilde x)\dvtx x,\tilde x\in\mathbb R^d\}$. Then the
functions\vadjust{\goodbreak}
$(x,\tilde x,\omega)\mapsto W(x,\tilde x)(\omega),\tilde W(x,\tilde
x)(\omega)$ defined by
\begin{eqnarray*}
W(x,\tilde x)&:=& %
\cases{Y(x,\tilde x),&\quad$\mbox{if }\mathcal U\leq
\alpha \bigl(x,Y(x,\tilde x)\bigr),$\vspace*{2pt}
\cr
x,&\quad$\mbox{if }\mathcal U>
\alpha\bigl(x,Y(x,\tilde x)\bigr)$,} %
\\
\tilde W(x,\tilde x)&:=& %
\cases{\tilde Y(x,\tilde x),&\quad$\mbox{if }
\mathcal U\leq\alpha \bigl(x,\tilde Y(x,\tilde x)\bigr),$\vspace*{2pt}
\cr
\tilde
x,&\quad$\mbox{if }\mathcal U>\alpha \bigl(x,\tilde Y(x,\tilde x)\bigr),$} %
\end{eqnarray*}
realize a Markovian coupling between the Metropolis--Hastings
transition functions $q(x,\cdot),x\in\mathbb R^d$, that is,
\[
W(x,\tilde x) \sim q(x,\cdot)\quad \mbox{and}\quad \tilde W(x,\tilde x) \sim q(\tilde x,
\cdot)
\]
for any $x,\tilde x\in\mathbb R^d$. This coupling is optimal in the
acceptance step in the sense that it minimizes the probability that a
proposed move from $x$ to $Y(x,\tilde x)$ is accepted and the
corresponding proposed move from $\tilde x$ to $\tilde Y(x,\tilde x)$
is rejected or vice versa.

%
%le2.5 #&#
\begin{lemma}[(Basic contractivity lemma for MH kernels)]\label{lemD}
$\!\!\!$For any \mbox{$x,\tilde x\in\mathbb R^d$},
\begin{eqnarray*}
&&\mathbb{E}\bigl[d\bigl(W(x,\tilde x),\tilde W(x,\tilde x)\bigr)\bigr]
\\
&&\qquad\leq \mathbb{E}\bigl[d\bigl(Y(x,\tilde x),\tilde Y(x,\tilde x)\bigr)\bigr]
\\
&&\qquad\quad{}+ \mathbb{E}\bigl[\bigl(d(x,\tilde x)-d\bigl(Y(x,\tilde x),\tilde Y(x,\tilde
x)\bigr)\bigr)\\
&&\hspace*{56pt}{}\times\max\bigl(1-\alpha \bigl(x,Y(x,\tilde x)\bigr),1-\alpha\bigl(
\tilde x,\tilde Y(x,\tilde x)\bigr)\bigr)\bigr]
\\
&&\qquad\quad{}+ \mathbb{E}\bigl[d(x,Y(x,\tilde x)\cdot\bigl(\alpha\bigl(x,Y(x,\tilde x)
\bigr)-\alpha \bigl(\tilde x,\tilde Y(x,\tilde x)\bigr)\bigr)^+\bigr]
\\
&&\qquad\quad{}+ \mathbb{E}\bigl[d(\tilde x,\tilde Y(x,\tilde x)\cdot\bigl(\alpha
\bigl(x,Y(x,\tilde x)\bigr)-\alpha\bigl(\tilde x,\tilde Y(x,\tilde x)\bigr)
\bigr)^-\bigr].
\end{eqnarray*}
\end{lemma}

\begin{pf}
By the definition of $W$ and by the triangle inequality, we obtain the estimate
\begin{eqnarray*}
&&\mathbb{E}\bigl[d\bigl(W(x,\tilde x),\tilde W(x,\tilde x)\bigr)\bigr]
\\
&&\qquad\leq \mathbb{E}\bigl[d\bigl(Y(x,\tilde x),\tilde Y(x,\tilde x)\bigr); \mathcal
U<\min \bigl(\alpha\bigl(x,Y(x,\tilde x)\bigr),\alpha\bigl(\tilde x,\tilde Y(x,
\tilde x)\bigr)\bigr)\bigr]
\\
&&\hspace*{10pt}\qquad\quad{}+d(x,\tilde x)\cdot\P\bigl[\mathcal U\ge\min\bigl(\alpha\bigl(x,Y(x,\tilde x)
\bigr),\alpha\bigl(\tilde x,\tilde Y(x,\tilde x)\bigr)\bigr)\bigr]
\\
&&\qquad\quad{}+\mathbb{E}\bigl[d\bigl(x,Y(x,\tilde x)\bigr); \alpha\bigl(\tilde x,\tilde
Y(x,\tilde x)\bigr)\leq\mathcal U<\alpha\bigl(x,Y(x,\tilde x)\bigr)\bigr]
\\
&&\qquad\quad{}+\mathbb{E}\bigl[d\bigl(\tilde x,\tilde Y(x,\tilde x)\bigr); \alpha\bigl( x,
Y(x,\tilde x)\bigr)\leq\mathcal U<\alpha\bigl(\tilde x,\tilde Y(x,\tilde x)\bigr)
\bigr]. %&\leq&\mathbb{E}[d(Y(x,\tilde x),\tilde Y(x,\tilde x))]\\
%&&+d(x,\tilde x)\cdot\mathbb{E}[\max(1-\alpha(x,Y(x,\tilde x)),1-
%&&+\Delta\cdot\mathbb{E}[|\alpha(x,Y(x,\tilde x))-\alpha(\tilde x,
\end{eqnarray*}
The assertion now follows by conditioning on $Y$ and $\tilde Y$.
\end{pf}

%
%re2.6 #&#
\begin{rmk} (1) Note that the upper bound in Lemma~\ref{lemD} is
close to
an equality. Indeed, the only estimate in the proof is the triangle
inequality that
has been applied to bound $d(x,\tilde Y)$ by $d(x,\tilde x)+d(\tilde x
,\tilde Y)$
and $d(\tilde x, Y)$ by $d(x,\tilde x)+d(x,Y)$.\vadjust{\goodbreak}

(2) For the couplings and distances considered in this paper, $d(Y,
\tilde Y)$
will always be deterministic. Therefore, the upper bound in the lemma
simplifies to
%
%
%e2.3 #&#
\begin{eqnarray}\label{eqlemD}
\quad &&\mathbb{E}\bigl[d(W,\tilde W)\bigr]
\nonumber
\\
&&\qquad\leq d(Y,\tilde Y)+ \bigl(d(x,\tilde x)-d(Y,\tilde Y)\bigr)\cdot\mathbb{E}
\bigl[\max\bigl(1-\alpha (x,Y),1-\alpha(\tilde x,\tilde Y)\bigr)\bigr]
\\
&&\qquad\quad{}+ \mathbb{E}\bigl[d(x,Y) \bigl(\alpha(x,Y)-\alpha(\tilde x,\tilde Y)\bigr)^+
+d(\tilde x,\tilde Y) \bigl(\alpha(x,Y)-\alpha(\tilde x,\tilde Y)\bigr)^-\bigr].\nonumber
\end{eqnarray}
Here $\mathbb{E} [\max(1-\alpha
(x,Y),1-\alpha(\tilde x,\tilde Y))]$ is the probability that at least
one of the proposals
is rejected.

(3) If the metric $d$ is bounded with diameter $\Delta$, then the last
two expectations
in the
upper bound in Lemma~\ref{lemD} can be estimated by $\Delta$ times the
probability
$\mathbb{E} [|\alpha
(x,Y)-\alpha(\tilde x,\tilde Y)|]$ that one of the proposals is
rejected and the other one is accepted. Alternatively (and usually more
efficiently), these terms can be estimated
by H\"older's inequality.
\end{rmk}

%s3 #&#
\section{Contractivity of the proposal step}
\label{secProposal}
\label{secMALA}

In this section we assume $V\in C^2(\mathbb R^d)$.
We study contractivity properties of the Metropolis--Hastings
proposals defined in (\ref{eq24a}) and (\ref{eq32}).

Note first that the
\emph{Ornstein--Uhlenbeck proposals} do not depend on $V$. For $h\in
(0,2)$, the contractivity condition
%
%
%e3.1 #&#
\begin{equation}
\bigl\llVert Y_h^{\mathrm{OU}}(x)-Y_h^{\mathrm{OU}}(
\tilde x)\bigr\rrVert = \bigl\llVert (1-h/2 ) (x-\tilde x)\bigr\rrVert = (1-h/2 )
\| x-\tilde x\| \label{eq25}
\end{equation}
holds pointwise for any $x,\tilde x\in\mathbb R^d$ w.r.t. an arbitrary
norm $\| \cdot\|$ on $\mathbb R^d$.

For the \emph{semi-implicit Euler proposals}
\[
Y_h(x) = x-\frac h2\nabla U(x)+\sqrt{h-\frac{h^2}{4}}Z, \qquad Z\sim
\gamma^d.
\]
Wasserstein contractivity does not necessarily hold. Close to optimal sufficient
conditions for contractivity w.r.t. the minus norm can be obtained in
a straightforward way by considering the derivative of $Y_h$ w.r.t. $x$.
%
%
%le3.1 #&#
\begin{lemma}\label{lemO}
Let $h\in(0,2)$, and let $C$ be a convex subset of $\mathbb R^d$. If
there exists a constant $\lambda\in(0,\infty)$ such that
%
%
%e3.2 #&#
\begin{equation}
\biggl\llVert \biggl(I_d-\frac h2\nabla^2U(x) \biggr)
\cdot\eta\biggr\rrVert _- \leq \lambda\|\eta\|_- \qquad\mbox{for any }\eta\in\mathbb
R^d,x\in C, \label{eq58}
\end{equation}
then
\[
\bigl\|Y_h(x)-Y_h(\tilde x)\bigr\|_- \leq\lambda\|x-\tilde x\|_-
\qquad\mbox{for any }x,\tilde x\in C.
\]
\end{lemma}

\begin{pf}
If (\ref{eq58}) holds, then
\[
\bigl\|\partial_\eta Y_h(x)\bigr\|_- = \biggl\llVert \eta-\frac h2
\nabla^2U(x)\cdot \eta \biggr\rrVert _- \leq\lambda\|\eta\|_-\vadjust{\goodbreak}
\]
for any $x\in C$ and $\eta\in\mathbb R^d$. Hence
\begin{eqnarray*}
\bigl\|Y_h(x)-Y_h(\tilde x)\bigr\|_-&=&\biggl\llVert \int
_0^1\frac
{d}{dt}Y_h
\bigl(tx+(1-t)\tilde x\bigr) \,dt\biggr\rrVert _-
\\
&\le&\int_0^1\bigl\llVert
\partial_{x-\tilde x}Y_h\bigl(tx+(1-t)\tilde x\bigr)\bigr\rrVert
_- dt
\\
&\leq&\lambda\|x-\tilde x\|_- \qquad\mbox{for }x,\tilde x\in C.
\end{eqnarray*}
\upqed\end{pf}

%
%re3.2 #&#
\begin{rmk}
(1) Note that condition (\ref{eq58}) requires a bound on $\nabla^2U$ in
both directions. This is in contrast to the continuous time case where
a lower bound by a strictly positive constant is sufficient to
guarantee contractivity of the derivative flow.

(2) Condition (\ref{eq58}) is equivalent to
%
%
%e3.3 #&#
\begin{equation}
\xi\cdot\eta-\frac h2\partial_{\xi\eta}^2U(x) \leq\lambda\|\xi
\|_+\| \eta\|_- \qquad\mbox{for any } x\in C,\xi,\eta\in\mathbb R^d.
\label{eq59a}
\end{equation}
%
%where $\|\xi\|_+=\sup\{\xi\cdot\eta: \eta\in\mathbb R^d\mbox{ with }\|
\end{rmk}

%For $x,y\in\mathbb R^d$ let

Recall that for $R\in(0,\infty]$,
%
%
%e3.4 #&#
\begin{equation}
\label{eq58a} M(R) = \sup\bigl\{\bigl\|\nabla^2U(x)\cdot\eta\bigr\|_-\dvtx \eta\in
B_1^-, x\in B_R^-\bigr\}. %L_2^-(x,y) = \sup\{\|\partial^2V(z)\|_{-,-}: z\in[x,y]\},
\end{equation}
%
%where $\|\partial^2 V(z)\|_{-,-} = \sup\{\|\partial^2V(z)\eta\|_-:
%$\partial^2V(z)$ as an operator from $(\mathbb R^d,\| \cdot\|_-)$ to
%$(\mathbb R^d,\| \cdot\|_-)$. Clearly,

Hence $M(R)$ bounds the second derivative of $U$ on $B_R^-$ in both directions,
whereas the constant $K$ in Assumption~\ref{assA2} is a strictly
positive lower
bound for the second derivative. We also define
%
%
%e3.5 #&#
\begin{equation}
\label{eq58b} N(R) = \sup\bigl\{\bigl\|\nabla^2V(x)\cdot\eta\bigr\|_-\dvtx \eta\in
B_1^-, x\in B_R^-\bigr\}. %L_2^-(x,y) = \sup\{\|\partial^2V(z)\|_{-,-}: z\in[x,y]\},
\end{equation}
Note that $M(R)\le1+N(R)$. As a consequence of Lemma~\ref{lemO} we obtain:

%
%pr3.3 #&#
\begin{prop}\label{propP}
For any $h\in(0,2)$ and $x,\tilde x\in B_R^-$,
%
%
%e3.6 #&#
\begin{equation}
\label{eq59} \bigl\|Y_{h}(x)-Y_h(\tilde x)\bigr\|_- \leq
\biggl(1-\frac{1-{N(R)}}2 h \biggr)\cdot\|x-\tilde x\|_-.
\end{equation}
Moreover, if Assumption~\ref{assA2} holds, then
%$\| \cdot\|_- = \langle
% \cdot, \cdot\rangle^{1/2}$ for an inner product $\langle
% \cdot, \cdot\rangle$ on $\mathbb R^d$, and there exists a
%constant $K\in\mathbb R$ such that
%then
%
%
%e3.7 #&#
\begin{equation}
\label{eq61} \bigl\|Y_h(x)-Y_h(\tilde x)\bigr\|_- \leq \biggl(1-
\frac K2 h+\frac{M(R)^2}{8} {h^2} \biggr)\cdot\| x-\tilde x\|_-.
\end{equation}
\end{prop}

\begin{pf}
Note that for $z\in[x,\tilde x]$ and $\eta\in\mathbb R$,
%
%
%e3.8 #&#
\begin{equation}
\biggl(I-\frac h2\nabla^2U(z)\biggr)\cdot\eta= \biggl(1-\frac h2
\biggr)\eta-\frac h2\nabla ^2V(z)\cdot\eta. \label{eq64}
\end{equation}
Therefore, by (\ref{eq58b}),
\[
\biggl\| \biggl(I-\frac h2\nabla^2U(z)\biggr)\cdot\eta\biggr\|_- \leq\biggl(1-
\frac h2\biggr)\|\eta\| _-+\frac h2N(R)\cdot\|\eta\|_-.
\]
Inequality (\ref{eq59}) now follows by Lemma~\ref{lemO}.\vadjust{\goodbreak}

Moreover, if Assumption~\ref{assA2} holds, then for $z\in[x,\tilde x]$
and $\eta\in\mathbb R^d$,
\begin{eqnarray*}
\biggl\| \biggl(I-\frac h2\nabla^2U(z)\biggr)\cdot\eta
\biggr\|_-^2& = &\|\eta\|_-^2-h\bigl\langle\eta,
\nabla^2U(z)\cdot\eta\bigr\rangle +\frac
{h^2}{4}\bigl\|
\nabla^2U(z)\cdot\eta\bigr\|^2_-
\\
&\leq& \bigl(1-Kh+M(R)^2h^2/4\bigr) \|\eta
\|_-^2 \\
&=& \bigl(1-Kh/2+M(R)^2h^2/8
\bigr)^2 \|\eta\|^2_-.
\end{eqnarray*}
Inequality (\ref{eq61}) again follows by Lemma~\ref{lemO}.
\end{pf}

%s4 #&#
\section{Upper bounds for rejection probabilities}
\label{secRP}

In this section we derive the upper bounds for the MH rejection
probabilities stated in
Proposition~\ref{propI}.
As a first step we prove the explicit formula for the MALA acceptance
probabilities
w.r.t. explicit and semi-implicit Euler proposals
stated in Proposition~\ref{propE}:

\begin{pf*}{Proof of Proposition~\ref{propE}}
For explicit Euler proposals with given step size $ h >0$,
\begin{eqnarray*}
&&\log\gamma^d(x)p_ h ^{\mathrm{Euler}}(x,y)\\
&&\qquad =
\frac12 |x|^2+\frac1{2 h }\biggl\llvert y- \biggl(1-\frac h 2
\biggr)x+\frac h 2\nabla V(x)\biggr\rrvert ^2+C
\\
&&\qquad=\frac{1}{2 h } \biggl( h |x|^2+|y-x|^2- h x\cdot
y+\frac14 h ^2|x|^2+ h (y-x)\cdot\nabla V(x)
\\
&&\hspace*{120pt}\qquad\quad{}+\frac12 h ^2 x\cdot\nabla V(x)+\frac14 h ^2\bigl|
\nabla V(x)\bigr|^2 \biggr)+C
\\
&&\qquad=S(x,y)+\frac12(y-x)\cdot\nabla V(x)+\frac h 8\bigl|x+\nabla V(x)\bigr|^2
\end{eqnarray*}
with a normalization constant $C$ that does not depend on $x$ and $y$,
and a symmetric function $S\dvtx \mathbb R^d\times\mathbb R^d\to\mathbb R$.
Therefore, by (\ref{eq21}),
\begin{eqnarray*}
G_ h ^{\mathrm{Euler}}(x,y)&=&V(y)-V(x)+\log\gamma^d(x)p_ h
^{\mathrm
{Euler}}(x,y)-\log\gamma^d(y)p_ h
^{\mathrm{Euler}}(y,x)
\\
&=&V(y)-V(x)-(y-x)\cdot\bigl(\nabla V(y)+\nabla V(x)\bigr)/2
\\
&&{}+ h \bigl(\bigl|y+\nabla V(y)\bigr|^2-\bigl|x+\nabla V(x)\bigr|^2\bigr)/8.
\end{eqnarray*}
Similarly, for semi-implicit Euler proposals we obtain
\begin{eqnarray*}
-\log\gamma^d(x)p_h(x,y) &=&
\frac12|x|^2+\frac 12\biggl\llvert y- \biggl(1-\frac h2 \biggr)x+\frac
h2\nabla V(x)\biggr\rrvert ^2 \Big/ \biggl(h-\frac{h^2}4 \biggr)+C
\\
&=&\frac12 \biggl( \biggl(h-\frac{h^2}4 \biggr)|x|^2+
\biggl\llvert y- \biggl(1-\frac h2 \biggr)x+\frac h2\nabla V(x)\biggr\rrvert
^2 \biggr) \\
&&{}\Big/ \biggl(h-\frac{h^2}4 \biggr)+C
\\
&=&-\frac12\frac{h}{4-h}|x|^2+\tilde S(x,y)+\frac12\cdot
\frac
{4}{4-h}(y-x)\nabla V(x)\\
&&{}+\frac12\cdot\frac{h}{4-h}\bigl|x+\nabla
V(x)\bigr|^2
\\
&=&\tilde S(x,y)+\frac12(y-x)\cdot\nabla V(x)\\
&&{}+\frac12\frac
{4}{4-h}
\bigl[(y+x)\cdot\nabla V(x)+\bigl|\nabla V(x)\bigr|^2 \bigr],
\end{eqnarray*}
and, therefore,
\begin{eqnarray*}
&&G_h(x,y)\\
&&\qquad=V(y)-V(x)-(y-x)\cdot\bigl(\nabla V(y)+\nabla V(x)
\bigr)/2
\\
&&\qquad\quad{}+\frac{h}{8-2h} \bigl[(y+x)\cdot\bigl(\nabla V(y)-\nabla V(x)\bigr)+\bigl|
\nabla V(y)\bigr|^2-\bigl|\nabla V(x)\bigr|^2 \bigr].\quad
\end{eqnarray*}
\upqed\end{pf*}

From now on we assume that Assumption~\ref{assA1} holds. We will derive
upper bounds for the functions
$G_h(x,y)$ computed in Proposition~\ref{propE}. By (\ref{eq8}), these
directly imply corresponding
upper bounds for the MALA rejection probabilities.

Let $\partial^n_{\xi_1,\ldots, \xi_n}V (z)$ denote the $n$th-order directional derivative of the function $V$ at~$z$
in directions ${\xi_1,\ldots, \xi_n}$. By $\partial^nV$ we
denote the $n$th-order differential of $V$, that is,
the $n$-form $(\xi_1,\ldots, \xi_n)\mapsto
\partial^n_{\xi_1,\ldots, \xi_n}V$.
%We fix a norm
%$\| \cdot\|_-$ on $\mathbb R^d$ such that
%We will assume regularity of $V$ w.r.t. $\| \cdot\|_-$. For $x,y\in
%denote the line segment between $x$ and $y$.
For $x,\tilde x\in\mathbb R^d$ and $n=1,2,3, 4$ let
%
%
%e4.1 #&#
\begin{equation}
L_n(x,\tilde x) = \sup\bigl\{\bigl(\partial^n_{\xi_1,\ldots, \xi_n}V
\bigr) (z)\dvtx z\in [x,\tilde x],\xi_1,\ldots,\xi_n\in
B_1^{-}\bigr\}. \label{eq38}
\end{equation}
In other words,
\[
L_n(x,\tilde x) = \sup_{z\in[x,\tilde x]}\bigl\|\bigl(
\partial^nV\bigr) (z)\bigr\|_-^\ast,
\]
where $\| \cdot\|_-^\ast$ is the dual norm on $n$-forms defined by
\[
\|l\|_-^\ast= \sup\bigl\{l(\xi_1,\ldots,
\xi_n)\dvtx \xi_1,\ldots,\xi _n\in
B_1^-\bigr\}.
\]
In particular,
\[
L_1(x,\tilde x) = \sup_{z\in[x,\tilde x]}\bigl\|\nabla V(z)\bigr\|_+.
\]
By Assumption~\ref{assA1},
%
%
%e4.2 #&#
\begin{equation}\qquad
L_n(x,\tilde x) \leq C_n\cdot\max\bigl(1,\|x\|_-,\|\tilde x
\|_-\bigr)^{p_n}\qquad \forall x,\tilde x\in\mathbb R^d,n=1,2,3,4.
\label{eq40}
\end{equation}

%In particular,
%where $\nabla$ denotes the Euclidean gradient, and
%A corresponding representation with $\nabla$ replaced by $\nabla^n$
%holds for $L_n(x,y)$. Note that if $\| \cdot\|_-$ is weaker than the
%Euclidean norm, then $\| \cdot\|_+$ is stronger than the Euclidean
%norm.

We now derive upper bounds for the terms in the expression for
$G_h(x,y)$ given in Proposition~\ref{propE}. We first express the
leading order term in terms of $3$rd derivatives of $V$:

%
%le4.1 #&#
\begin{lemma}\label{lemF}
For $x,y\in\mathbb R^d$,
\begin{eqnarray*}
&& V(y)-V(x)-\frac{y-x}{2}\cdot\bigl(\nabla V(y)+\nabla V(x)\bigr) \\
&&\qquad= -\frac12
\int_0^1t(1-t)\partial^3_{y-x}V
\bigl((1-t)x+ty\bigr) \,dt.
\end{eqnarray*}
\end{lemma}

\begin{pf}
A second-order expansion for $f(t)=V(x+t(y-x)), t\in[0,1]$, yields
\begin{eqnarray*}
V(y)-V(x)&=&\int_0^1\partial_{y-x}V
\bigl(x+t(y-x)\bigr) \,dt
\\
&=&(y-x)\cdot\nabla V(x)+\int_0^1\int
_0^t\partial^2_{y-x}V
\bigl(x+s(y-x)\bigr) \,ds \,dt
\\
&=&(y-x)\cdot\nabla V(x)+\int_0^1(1-s)
\partial^2_{y-x}V\bigl(x+s(y-x)\bigr) \,ds,
\end{eqnarray*}
and, similarly,
\begin{eqnarray*}
V(y)-V(x)&=&(y-x)\cdot\nabla V(y)-\int_0^1
\int_t^1\partial ^2_{y-x}V
\bigl(x+s(y-x)\bigr) \,ds \,dt
\\
&=&(y-x)\cdot\nabla V(y)-\int_0^1s
\partial^2_{y-x}V\bigl(x+s(y-x)\bigr) \,ds.
\end{eqnarray*}
By averaging both equations, we obtain
\begin{eqnarray*}
&&V(y)-V(x)-\frac{y-x}{2}\cdot\bigl(\nabla V(y)+\nabla V(x)\bigr)
\\
&&\qquad=\frac12\int_0^1(1-2s)\partial^2_{y-x}V
\bigl(x+s(y-x)\bigr) \,ds
\\
&&\qquad=\frac12 \int_0^1t(1-t)
\partial^3_{y-x}V\bigl(x+t(y-x)\bigr) \,dt.
\end{eqnarray*}
Here we have used that for any function $g\in C^1([0,1])$,
\begin{eqnarray*}
\int_0^1(1-2s)g(s) \,ds&=&\int
_0^1(1-2s)\int_0^sg'(t)
\,dt \,ds
\\
&=&\int_0^1\int_t^1(1-2s)
\,ds g'(t) \,dt = -\int_0^1t(1-t)g'(t)
\,dt.
\end{eqnarray*}
\upqed\end{pf}

%
%le4.2 #&#
\begin{lemma}\label{lemG}
For $x,y\in\mathbb R^d$, the following estimates hold:
\begin{longlist}[(1)]
\item[(1)] $|V(y)-V(x)|\le L_1(x,y)\cdot\|y-x\|_-$;
\item[(2)] $\llvert V(y)-V(x)-\frac{y-x}{2}\cdot(\nabla V(y)+\nabla
V(x))\rrvert  \leq\frac1{12}L_3(x,y)\cdot\|y-x\|_-^3$;
\item[(3)] $|(\nabla U(y)+\nabla U(x))\cdot(\nabla V(y)-\nabla
V(x))|$
$ \leq L_2(x,y)\cdot\|\nabla U(y)+\nabla U(x)\|_-\cdot\|y-x\|_-$;
\item[(4)] $\|\nabla U(y)+\nabla U(x)\|_- \leq2\|\nabla U(x)\|
_-+(1+L_2(x,y))\cdot\|y-x\|_-$.
\end{longlist}
\end{lemma}

%
%re4.3 #&#
\begin{rmk}\label{rmexplicit}
The estimates in Lemma~\ref{lemG} provide a bound for the terms in the
expression (\ref{eq33}) for $G_h(x,y)$ in the case of semi-implicit
Euler proposals. For explicit Euler proposals, one also has to bound
the term
\[
\bigl|\nabla U(y)\bigr|^2-\bigl|\nabla U(x)\bigr|^2 = \bigl|y+\nabla
V(y)\bigr|^2-\bigl|x+\nabla V(x)\bigr|^2.
\]
Note that even when $V$ vanishes, this term cannot be controlled in
terms of $\| \cdot\|_-$ in general. A valid upper bound is
\[
\bigl|\nabla U(y)+\nabla U(x)\bigr|\cdot|y-x|+L_2(x,y)\bigl\|\nabla U(y)+\nabla
U(x)\bigr\| _-\cdot\|y-x\|_-.
\]
\end{rmk}

\begin{pf*}{{Proof of Lemma~\ref{lemG}}}
By Lemma~\ref{lemF} and by definition of $L_n(x,y)$, we obtain
\begin{eqnarray*}
&&\bigl|V(y)-V(x)\bigr| \leq \sup_{z\in[x,y]}\bigl|\partial_{y-x}V(z)\bigr|
\le L_1(x,y)\cdot\|y-x\|_-,
\\
&&\biggl\llvert V(y)-V(x)-\frac{y-x}{2}\cdot\bigl(\nabla V(y)+\nabla V(x)
\bigr)\biggr\rrvert
\\
&&\qquad\leq \frac12\int_0^1t(1-t) \,dt \sup
_{z\in[x,y]}\bigl|\partial^3_{y-x}V(z)\bigr| \leq
\frac1{12}L_3(x,y)\cdot\|y-x\|_-^3,\\
&&\bigl(\nabla U(y)+\nabla U(x)\bigr)\cdot\bigl(\nabla V(y)-\nabla V(x)
\bigr)
\\
&&\qquad=\partial_{\nabla U(y)+\nabla U(x)}V(y)-\partial_{\nabla
U(y)+\nabla
U(x)}V(x)
\\
&&\qquad= \int_0^1\partial^2_{\nabla U(y)+\nabla U(x),y-x}V
\bigl((1-t)x+ty\bigr) \,dt
\\
&&\qquad\leq L_2(x,y)\cdot\bigl\|\nabla U(y)+\nabla U(x)\bigr\|_-\cdot\|y-x\|_-.
\end{eqnarray*}
Moreover, the estimate
\begin{eqnarray*}
\bigl\|\nabla U(y)+\nabla U(x)\bigr\|_-&\leq&2\bigl\|\nabla U(x)\bigr\|_-+\bigl\|\nabla U(y)-\nabla U(x)
\bigr\|_-
\\
&\leq&2\bigl\|\nabla U(x)\bigr\|_-+\|y-x\|_-+\bigl\|\nabla V(y)-\nabla V(x)\bigr\|_-
\\
&\leq&2\bigl\|\nabla U(x)\bigr\|_-+\bigl(1+L_2(x,y)\bigr)\cdot\|y-x\|_-
\end{eqnarray*}
holds by definition of $L_2(x,y)$ and since
\begin{eqnarray*}
\hspace*{20pt}\bigl\|\nabla V(y)-\nabla V(x)\bigr\|_-&\leq&\bigl|\nabla V(y)-\nabla V(x)\bigr| = \sup
_{|\xi|=1}\bigl(\partial_\xi V(y)-\partial_\xi
V(x)\bigr)
\\
&\leq&\sup_{\|\xi\|_-\le1}\bigl(\partial_\xi V(y)-
\partial_\xi V(x)\bigr).\hspace*{100pt}\qed
\end{eqnarray*}
\noqed\end{pf*}
%We now prove an upper bound for the average rejection probability in
%the next step of semi-implicit MALA when the chain is currently at
%position $x\in\mathbb R^d$.
%For this purpose, we assume from now on that $V$ is sufficiently
%smooth with derivatives growing at most polynomially w.r.t. $\| \cdot
%
%The function $V$ is in $C^4(\mathbb R^d)$, and there exist finite
%constants $C_n,p_n\in[0,\infty)$ such that
%for any $x\in\mathbb R^d$, $\xi_1,\ldots,\xi_n\in S_-^{d-1}$ and
%$n=2,3,4$.
Recalling the definitions of $Y_h^{\mathrm{OU}}(x)$ and $Y_h(x)$ from
(\ref{eq24a}) and (\ref{eq32}), we obtain:

%
%le4.4 #&#
\begin{lemma}\label{lemH}
For $x\in\mathbb R^d,h\in(0,2)$ and $n\in\{ 1,2,3,4\}$ with $p_n\ge1$,
we have:
\begin{longlist}[(1)]
\item[(1)] $\|Y_h^{\mathrm{OU}}(x)-x\|_- \leq\frac h2\| x\|_-+\sqrt
h\|
Z\|_-$;
\item[(2)] $\|Y_h(x)-x\|_- \leq\frac h2\|\nabla U(x)\|_-+\sqrt h\|Z\|_-$;
\item[(3)] $\|Y_h^{\mathrm{OU}}(x)\|_- \leq(1-\frac h2)\| x\|
_-+\sqrt
h\|Z\|_-$;
\item[(4)] $\|Y_h(x)\|_- \leq\|x\|_-+\frac h2\|\nabla U(x)\|_-+\sqrt
h\|Z\|_-$;
\item[(5)] $L_n(x,Y_h^{\mathrm{OU}}(x)) \leq C_n2^{p_n-1} (\max
(1,\|
x\|_-)^{p_n}+h^{p_n/2}\|Z\|^{p_n}_- )$;
\item[(6)] $L_n(x,Y_h(x)) \leq C_n3^{p_n-1} (\max(1,\|x\|
_-)^{p_n}+ (\frac h2 )^{p_n}\|\nabla U(x)\|
_-^{p_n}+\break   h^{p_n/2}\|
Z\|^{p_n}_- )$.
\end{longlist}
\end{lemma}

\begin{pf}
Estimates (1)--(4) are direct consequences of the triangle inequality.
Moreover, by (3) and (4),
%and hence
%
\[
\max\bigl(1,\|x\|_-,\bigl\|Y_h^{\mathrm{OU}}(x)\bigr\|_-\bigr) \leq\max\bigl(1,
\|x\| _-\bigr)+\sqrt h\|Z\|_-
\]
and
\[
\max\bigl(1,\|x\|_-,\bigl\|Y_h(x)\bigr\|_-\bigr) \leq\max\bigl(1,\|x\|_-\bigr)+\frac h2
\bigl\|\nabla U(x)\bigr\|_-+\sqrt h\|Z\|_-.
\]
Estimates (5) and (6) now follow from (\ref{eq40}) and H\"older's inequality.
\end{pf}

We now combine the estimates in Lemmas \ref{lemG} and~\ref{lemH}
in order
to prove Proposition~\ref{propI} and the first part of Proposition~\ref
{propMOU}:

\begin{pf*}{Proof of Proposition~\ref{propI}}
By (\ref{eq8}) and Proposition~\ref{propE}, for \mbox{$h\in(0,2)$},
%
%
%e4.3 #&#
\begin{equation}
\mathbb{E} \bigl[\bigl(1-\alpha_h\bigl(x,Y_h(x)\bigr)
\bigr)^k \bigr]^{1/k} \le\bigl\llVert G_h
\bigl(x,Y_h(x)\bigr)^+\bigr\rrVert _{L^k} \leq I+\frac h4\mathit{II},
\label{eq41}
\end{equation}
where
\begin{eqnarray*}
I&=&\biggl\llVert V\bigl(Y_h(x)\bigr)-V(x)-\frac{Y_h(x)-x}2\cdot
\bigl(\nabla V\bigl(Y_h(x)\bigr)-\nabla V(x)\bigr)\biggr\rrVert
_{L^k},
\\
\mathit{II}&=&\bigl\llVert \bigl(\nabla U\bigl(Y_h(x)\bigr)+\nabla U(x)
\bigr)\cdot\bigl(\nabla V\bigl(Y_h(x)\bigr)-\nabla V(x)\bigr)\bigr
\rrVert _{L^k}.
\end{eqnarray*}
By Lemma~\ref{lemG},
%
%e4.4 #&#
\begin{eqnarray}\quad
\label{eq42} I& \leq& \mathbb{E} \bigl[L_3\bigl(x,Y_h(x)
\bigr)^k\cdot\|Y_h(x)-x\|_-^{3k}
\bigr]^{1/k}/12 \quad\mbox{and }
\\
\label{eq43} \mathit{II}&\leq&\mathbb{E} \bigl[L_2\bigl(x,Y_h(x)
\bigr)^k\cdot\|Y_h(x)-x\|_-^k
\nonumber
\\[-8pt]
\\[-8pt]
\nonumber
&&\hspace*{10pt}{} \times\bigl(2\bigl\|\nabla U(x)\bigr\|_-+\bigl(1+L_2
\bigl(x,Y_h(x)\bigr)\bigr)\cdot\|Y_h(x)-x\| _-
\bigr)^k \bigr]^{1/k}.
\end{eqnarray}
The assertion of Proposition~\ref{propI} for semi-implicit Euler
proposals is now a direct consequence of the estimates (2) and (6) in
Lemma~\ref{lemH}.\vadjust{\goodbreak} The assertion for
Ornstein--Uhlenbeck proposals follows similarly from (\ref{eq23}), the
estimates
(1) in Lemma~\ref{lemG}
and (1), (3) and (5) in Lemma~\ref{lemH}.
\end{pf*}

It is possible to write down the polynomial in Proposition~\ref{propI}
explicitly. For
semi-implicit Euler proposals, we illustrate this in the case $k=1$ and
$p_2=p_3=0$. Here, by (\ref{eq42}) and (\ref{eq43}) we obtain
\begin{eqnarray*}
I&\leq&\frac{C_3}{12}\mathbb{E} \bigl[ \bigl( h\bigl\|\nabla U(x)\bigr\| _-/2+\sqrt h
\|Z\|_- \bigr)^3 \bigr] \leq\frac{C_3}{4} \bigl(
h^3\bigl\|\nabla U(x)\bigr\| ^3_-/8+h^{3/2}m_3
\bigr),
\\
\mathit{II}&\leq&C_2 \mathbb{E} \bigl[ \bigl( h\bigl\|\nabla U(x)\bigr\|_- /2+\sqrt h\|
Z\| _- \bigr)
\\
&&\hspace*{20pt}{} \times \bigl(2\bigl\|\nabla U(x)\bigr\|_-+(1+C_2) \bigl( h\bigl\|\nabla
U(x)\bigr\| _-/2+\sqrt h\|Z\|_- \bigr) \bigr) \bigr]
\\
&\leq&C_2 \biggl(h\bigl\|\nabla U(x)\bigr\|_-^2+2\sqrt h\bigl\|\nabla
U(x)\bigr\| _-m_1\\
&&\hspace*{17pt}{}+(1+C_2) \biggl(\frac{h^2}2\bigl\|\nabla
U(x)\bigr\|^2_-+2hm_2\biggr) \biggr).
\end{eqnarray*}
Hence by (\ref{eq41}),
%
%
%e4.6 #&#
\begin{eqnarray}\label{eqPE}
\mathbb{E} \bigl[1-\alpha_h\bigl(x,Y_h(x)\bigr)\bigr]&
\leq&h^{3/2}\cdot \biggl(\frac 14C_3m_3+
\frac12C_2m_1\bigl\|\nabla U(x)\bigr\|_- \biggr)
\nonumber
\\
&&{}+h^2\cdot \biggl(\frac14 C_2\bigl\|\nabla U(x)
\bigr\|^2_-+\frac 12C_2(1+C_2)m_2
\biggr)
\\
&&{}+h^3\cdot \biggl(\frac1{32}C_3\bigl\|\nabla U(x)
\bigr\|^3_-+\frac 18C_2(1+C_2\bigl)\bigl\| \nabla U(x)
\bigr\|^2_- \biggr).\nonumber
\end{eqnarray}

For Ornstein--Uhlenbeck proposals, we derive the explicit bound for the
rejection probabilities stated in Proposition~\ref{propMOU} for the case
$k=1$ and $p_2=0$.
\begin{pf*}{Proof of Proposition~\ref{propMOU}, first part}
If $p_2=0$, then for any $x\in\mathbb R^d$,
%
%
%e4.7 #&#
\begin{equation}\quad
\bigl\|\nabla V (x)\bigr\|_+ \le\bigl \|\nabla V (0)\bigr\|_+ +\bigl\|\nabla V (x)-\nabla V (0)\bigr\|_+ \le
C_1+C_2\cdot\| x\|_-.\label{eqnablav}
\end{equation}
Therefore, for any $x,y\in\mathbb R^d$,
\[
\bigl|V(y)-V(x)\bigr| \le \bigl( C_1+C_2\cdot\max\bigl(\| x\|_-,\| y
\|_-\bigr) \bigr)\cdot\| y-x\|_-.
\]
Hence, by (\ref{eq23}) and by (1) and (3) in Lemma~\ref{lemH},
\begin{eqnarray*}
&&\mathbb{E} \bigl[1-\alpha^\mathrm{OU}\bigl(x,Y_h^\mathrm
{OU}(x)\bigr) \bigr] \\
&&\qquad\le \mathbb{E} \bigl[\bigl(V\bigl(Y_h^\mathrm
{OU}(x)\bigr)-V(x)\bigr)^+ \bigr]
\\
&&\qquad\le \mathbb{E} \bigl[ \bigl( C_1+C_2\cdot\max\bigl(\|
x\|_-,\bigl\| Y_h^\mathrm {OU}(x)\bigr\|_-\bigr) \bigr)\cdot \bigl\|
Y_h^\mathrm{OU}(x)-x\bigr\|_- \bigr]
\\
&&\qquad\le \mathbb{E} \bigl[ \bigl( C_1+C_2\cdot\bigl(\| x\|_-+
\sqrt h\| Z\| _-\bigr) \bigr)\cdot \bigl(h\| x\|_-/2+\sqrt h\| Z\|_- \bigr) \bigr]
\\
&&\qquad= m_1 \bigl(C_1+C_2 \| x\|_-\bigr)\cdot
h^{1/2}
\\
& &\qquad\quad{}+ \frac12 \bigl(2m_2C_2+C_1\| x
\|_-+C_2\| x\|_-^2\bigr)\cdot h + \frac12m_1C_2
\| x\|_-\cdot h^{3/2}.
\end{eqnarray*}
\upqed\end{pf*}

%s5 #&#
\section{Dependence of rejection on the current state}
\label{secDRP}
We now derive estimates for the derivatives of the functions
%
%
%e5.1 #&#
%e5.2 #&#
\begin{eqnarray}
F_h(x,w)&=&G_h \biggl(x, x-\frac h2\nabla U(x)+w
\biggr), \label{eq44}
\\
F_h^\mathrm{OU}(x,w)&=&G^\mathrm{OU}\biggl(x, x-\frac
h2 x+w\biggr), \qquad(x,w)\in \mathbb R^d\times\mathbb R^d,
\label{eq44OU}
\end{eqnarray}
w.r.t. $x$. Since
%
%
%e5.3 #&#
%e5.4 #&#
\begin{eqnarray}
G_h\bigl(x,Y_h(x)\bigr)&=&F_h\biggl(x,
\sqrt{h-\frac{h^2}4}Z\biggr)\qquad \mbox{ with }Z\sim \gamma ^d,
\mbox{ and} \label{eq45}
\\
G^\mathrm{OU}\bigl(x,Y_h^\mathrm{OU}(x)
\bigr)&=&F_h^\mathrm{OU}\biggl(x,\sqrt {h-\frac
{h^2}4}Z
\biggr) \qquad\mbox{ with }Z\sim\gamma^d, \label{eq45OU}
\end{eqnarray}
these estimates can then be applied to control the dependence of
rejection on the current state $x$.

For \emph{Ornstein--Uhlenbeck proposals}, by (\ref{eq23}), we
immediately obtain
%
%
%e5.5 #&#
\begin{equation}
\label{eqKOU} \nabla_xF_h^\mathrm{OU}(x,w) =
\biggl(1-\frac h2\biggr) \bigl(\nabla V(y)-\nabla V(x) \bigr) -\frac h2\nabla
V(x),
\end{equation}
where $y:= ( 1-\frac h2) x+w$.

For \emph{semi-implicit Euler proposals}, the formula for the derivative
is more involved. To simplify the notation we set for $x\in\mathbb R^d$
and fixed $h\in(0,2)$,
%
%
%e5.6 #&#
\begin{equation}
x':= x-\frac h2\nabla U(x). \label{eq46}
\end{equation}
In the sequel, we use the conventions
\begin{eqnarray*}
v\cdot w&=&\sum_{i=1}^dv_iw_i,\qquad
(v\cdot T)_j=\sum_{i=1}^dv_iT_{i,j},\qquad
(T\cdot v)_j=\sum_{j=1}^dT_{i,j}v_j,\\
(S\cdot T)_{ik}&=&\sum_{j=1}^dS_{i,j}T_{j,k}
\end{eqnarray*}
for vectors $v,w\in\mathbb R^d$ and $(2,0)$ tensors $S,T\in\mathbb
R^d\otimes\mathbb R^d$. In particular,
\[
v\cdot(S\cdot T) = (v\cdot S)\cdot T,
\]
that is, the brackets may be omitted. We now give an explicit
expression for the derivative of $F_h(x,w)$ w.r.t. the first
variable:\vadjust{\goodbreak}

%
%pr5.1 #&#
\begin{prop}\label{propK}
Suppose $V\in C^2(\mathbb R^d)$. Then for any $x,w\in\mathbb R^d$,
\begin{eqnarray*}
\nabla_xF_h(x,w) &=& \nabla V(y)-\nabla V(x)-
\frac
{y-x}{2}\cdot \bigl(\nabla^2V(y)+\nabla^2V(x)
\bigr)
\\
&&{}+\frac h4(y-x)\cdot\nabla^2V(y)\cdot\bigl(I_d+
\nabla^2V(x)\bigr)
\\
&&{}+\frac h{8-2h}\bigl(\nabla V(y)-\nabla V(x)+\nabla U(y)+\nabla
U(x)\bigr)\\
&&{}\times
\bigl(\nabla^2V(y)-\nabla^2V(x)\bigr)
\\
&&{}-\frac{h^2}{16-4h}\bigl(\nabla V(y)-\nabla V(x)+\nabla U(y)+\nabla U(x)\bigr)
\cdot\nabla^2V(y)\\
&&{}\times \bigl(I_d+\nabla^2V(x)
\bigr)
\end{eqnarray*}
with $y:=x'+w=x-\frac h2\nabla U(x)+w$.
\end{prop}

\begin{pf}
Let
\[
W(x):= \nabla V(x) = \nabla U(x)-x,\qquad x\in\mathbb R^d.
\]
By Proposition~\ref{propE},
%
%
%e5.7 #&#
\begin{equation}\label{eq48}
F_h(x,w) = A_h(x,w)+\frac{h}{8-2h}B_h(x,w)
\end{equation}
for any $x,w\in \mathbb R^d$, where
\begin{eqnarray*}
A_h(x,w)&:=&V\bigl(x'+w\bigr)-V(x)-\frac{x'+w-x}{2}
\cdot\bigl(\nabla V\bigl(x'+w\bigr)+\nabla V(x)\bigr)\quad \mbox{and}
\\
B_h(x,w)&:=&\bigl(\nabla U\bigl(x'+w\bigr)+\nabla
U(x)\bigr)\cdot\bigl(\nabla V\bigl(x'+w\bigr)-\nabla V(x)\bigr).
\end{eqnarray*}
Noting that by (\ref{eq46}),
%
%
%e5.8 #&#
%e5.9 #&#
%e5.10 #&#
\begin{eqnarray}
\label{eq49}x-x'&=&\frac h2\nabla U(x) = \frac h2 x+\frac h2
\nabla V(x),
\\
\label{eq50}\nabla_x\bigl(x-x'\bigr)&=&\frac h2
\nabla^2U(x) = \frac h2 I_d+\frac h2\nabla^2V(x)\quad
\mbox{and}
\\
\label{eq51}\nabla_x x'&=&I_d-\frac h2
\nabla^2U(x) = \biggl(1-\frac h2 \biggr)I_d-\frac h2
\nabla^2V(x),
\end{eqnarray}
we obtain with $y=x'+w$
\begin{eqnarray*}
\nabla_xA_h(x,w)& =& W\bigl(x'+w\bigr)
\cdot\biggl(I_d-\frac h2\nabla^2U(x)\biggr)-W(x)
\\
&&{}-\frac{x'+w-x}{2}\cdot \biggl(\nabla W\bigl(x'+w\bigr)\cdot
\biggl(I_d-\frac h2\nabla ^2U(x)\biggr)+\nabla W(x)
\biggr)
\\
&&{}+\frac h4\bigl(W\bigl(x'+w\bigr)+W(x)\bigr)\cdot
\nabla^2U(x)
\\
&=&W(y)-W(x)-\frac{y-x}{2}\cdot\bigl(\nabla W(y)+\nabla W(x)\bigr)
\\
&&{}-\frac h4\bigl(W(y)-W(x)\bigr)\cdot\nabla^2U(x)+\frac h4(y-x)\cdot
\bigl(\nabla W(y)\cdot\nabla^2U(x)\bigr),\\ %\vspace*{-9pt}
\nabla_xB_h(x,w) &=& \bigl(W
\bigl(x'+w\bigr)-W(x)\bigr)\\
&&{}\times \biggl(\nabla ^2U
\bigl(x'+w\bigr)\cdot\biggl(I_d-\frac h2
\nabla^2U(x)\biggr)+\nabla^2U(x) \biggr)
\\
&&{}+\bigl(\nabla U\bigl(x'+w\bigr)+\nabla U(x)\bigr)\\
&&{}\times \biggl(
\nabla W\bigl(x'+w\bigr)\cdot \biggl(I_d-\frac h2
\nabla^2U(x)\biggr)-\nabla W(x) \biggr)
\\
&=&\bigl(W(y)-W(x)\bigr)\cdot\bigl(\nabla^2 U(y)\\
&&{}+
\nabla^2U(x)\bigr)+\bigl(\nabla U(y)+\nabla U(x)\bigr)\cdot\bigl(
\nabla W(y)-\nabla W(x)\bigr)
\\
&&{}-\frac h2\bigl(W(y)-W(x)\bigr)\cdot\bigl(\nabla^2U(y)\cdot
\nabla^2U(x)\bigr)\\
&&{}-\frac h2\bigl(\nabla U(y)+\nabla U(x)\bigr)\cdot
\bigl(\nabla W(y)\cdot\nabla^2U(x)\bigr).
\end{eqnarray*}
In total, we obtain
\begin{eqnarray*}
\nabla_xF_h(x,w) &=&\nabla_xA_h(x,w)+
\frac{h}{8-2h}\nabla _xB_h(x,w)
\\
&=&W(y)-W(x)-\frac{y-x}2\cdot\bigl(\nabla W(y)+\nabla W(x)\bigr)
\\
&&{}+\frac{h}{8-2h}\bigl(W(y)-W(x)\bigr)\cdot\bigl(\nabla^2U(y)-
\nabla^2U(x)\bigr)\\
&&{}+\frac h4(y-x)\cdot\nabla W(y)\cdot
\nabla^2U(x)
\\
&&{}+\frac{h}{8-2h}\bigl(\nabla U(y)+\nabla U(x)\bigr)\cdot\bigl(\nabla W(y)-
\nabla W(x)\bigr)
\\
&&{}+ \biggl(\frac{2h}{8-2h}-\frac h4 \biggr) \bigl(W(y)-W(x)\bigr)\cdot\nabla
^2U(x)
\\
&&{}-\frac{h^2}{16-4h}\bigl(W(y)-W(x)\bigr)\cdot\nabla^2U(y)\cdot
\nabla^2U(x)
\\
&&{}-\frac{h^2}{16-4h}\bigl(\nabla U(y)+\nabla U(x)\bigr)\cdot\nabla W(y)\cdot
\nabla^2U(x)
\\
&=&W(y)-W(x)-\frac{y-x}{2}\cdot\bigl(\nabla W(y)+\nabla W(x)\bigr)\\
&&{}+\frac
h4(y-x)\cdot\nabla W(y)\cdot\nabla^2U(x)
\\
&&{}+\frac h{8-2h}\bigl(W(y)-W(x)+\nabla U(y)+\nabla U(x)\bigr)\\
&&{}\times\bigl(\nabla
W(y)-\nabla W(x)\bigr)
\\
&&{}-\frac{h^2}{16-4h}\bigl(\bigl(W(y)-W(x)\bigr)\cdot\bigl(\nabla^2U(y)-I_d
\bigr)\\
&&{}+\bigl(\nabla U(y)+\nabla U(x)\bigr)\cdot\nabla W(y)\bigr)\cdot
\nabla^2U(x).
\end{eqnarray*}
Here we have used that
\[
\nabla^2U = I_d+\nabla^2V = I_d+
\nabla W.
\]
The assertion follows by applying this identity to the remaining
$\nabla
^2U$ terms as well.
\end{pf}

Similar to Lemma~\ref{lemG} above, we now derive bounds for the
individual summands in the expressions
for $\nabla_xF_h^\mathrm{OU}$ and $\nabla_xF_h$ in (\ref{eqKOU}) and
Proposition~\ref{propK}.
%Recall that
%denotes the dual norm of $\| \cdot\|_-$ on $\mathbb R^d$. Hence
%Note also that for a function $F\in C^1(\mathbb R^d)$,
%|F(y)-F(x)|&=&\left|\int_0^1(y-x)\cdot\nabla F((1-t)x+ty) dt\right|\\
%&\leq&\|y-x\|_-\cdot\sup_{z\in[x,y]}\|\nabla F(z)\|_+,
%i.e., the plus norm of $\nabla F$ determines the Lipschitz constant
%w.r.t. the minus norm.

%
%le5.2 #&#
\begin{lemma}\label{lemL}
For $V\in C^4(\mathbb R^d)$ and $x,y\in\mathbb R^d$ the following
estimates hold:
\begin{longlist}[(1)]
\item[(1)] $\llVert \nabla V(y)-\nabla V(x)\rrVert _+ \leq
L_2(x,y)\cdot
\|y-x\|_-$;
\item[(2)] $\llVert \nabla V(y)-\nabla V(x)-\frac{y-x}2\cdot(\nabla
^2V(y)+\nabla^2V(x))\rrVert _+ \leq
L_4(x,y)\cdot\|y-x\|_-^3/12$;
\item[(3)] $\|(y-x)\cdot\nabla^2V(y)\cdot(I_d+\nabla^2V(x))\|_+
\leq L_2(y,y)\cdot(1+L_2(x,x))\cdot\|y-x\|_-$;

\item[(4)] $\|(\nabla V(y)-\nabla V(x)+\nabla U(y)+\nabla U(x))\cdot
(\nabla^2V(y)-\nabla^2V(x))\|_+$
$ \leq L_3(x,y)\cdot(L_2(x,y) \|y-x\|_-+\|\nabla U(y)+\nabla U(x)\|
_-)\cdot\|y-x\|_-$;
\item[(5)] $\|(\nabla V(y)-\nabla V(x)+\nabla U(y)+\nabla U(x))\cdot
\nabla^2V(y)\cdot(I_d+\nabla^2V(x))\|_-$
$ \leq L_2(y,y)\cdot(L_2(x,y) \|y-x\|_-+\|\nabla U(y)+\nabla U(x)\|_-)
\cdot(1+L_2(x,x))$.
\end{longlist}
\end{lemma}

\begin{pf}
(1) For any $\xi\in\mathbb R^d$, we have
\[
\bigl|\partial_\xi V(y)-\partial_\xi V(x)\bigr| \le \sup
_{z\in
[x,y]}\bigl|\partial ^2_{y-x,\xi} V(z)\bigr| \le
L_2(x,y)\| x-y\|_-\|\xi\|_-.
\]
This proves (1) by definition of $\|\cdot\|_+$.

(2) By Lemma~\ref{lemF} applied to $\partial_\xi V$,
\begin{eqnarray*}
&&\biggl|\partial_\xi V(y)-\partial_\xi V(x)-
\frac{y-x}2\cdot \bigl(\nabla \partial_\xi V(y)- \nabla
\partial_\xi V(x)\bigr)\biggr|
\\
&&\qquad\le \frac12 \int_0^1 t(1-t) \,dt\cdot\sup
_{z\in[x,y]}\bigl|\partial ^3_{y-x}
\partial_\xi V(z)\bigr| \le \frac1{12}L_4(x,y)\| x-y
\|_-^3\|\xi \| _-.
\end{eqnarray*}

(3) For $v,w\in\mathbb R^d$, we have
%
%
%e5.11 #&#
\begin{equation}
\label{eqxxx} \bigl|v\cdot\nabla^2V(y)w\bigr| = \bigl|\partial^2_{vw}V(y)\bigr|
\le L_2(x,y)\| v\| _-\| w\|_-.
\end{equation}
Since $\|\cdot\|_-$ is weaker than the Euclidean norm, we obtain
\begin{eqnarray*}
&&\bigl|(y-x)\cdot\nabla^2V(y)\cdot\bigl(I+\nabla^2V(x)\bigr)
\cdot\xi\bigr| \\
&&\qquad\le L_2(y,y)\| y-x\|_-\bigl\| \bigl(I+\nabla^2V(x)
\bigr)\cdot\xi\bigr\|_-
\\
&&\qquad\le L_2(y,y)\| y-x\|_- \bigl(1+L_2(x,x)\bigr)\| \xi
\|_-.
\end{eqnarray*}

(4), (5) For $v,w\in\mathbb R^d$,
\begin{eqnarray*}
\bigl|v\cdot\bigl(\nabla^2V(y)-\nabla^2V(x)\bigr)\cdot w\bigr| &=&
\biggl\llvert \int_0^1\partial^3_{y-x,v,w }V
\bigl((1-t)x+ty\bigr) \,dt\biggr\rrvert
\\
&\le& L_3(x,y)\| y-x\|_-\| v\|_-\| w\|_-.
\end{eqnarray*}
Therefore,
\begin{eqnarray*}
\hspace*{-4pt}&&\bigl|\bigl(\nabla V(y)-\nabla V(x)+\nabla U(y)+\nabla U(x)\bigr)\cdot
\bigl(\nabla ^2V(y)-\nabla^2V(x)\bigr)\cdot \xi\bigr|
\\
\hspace*{-4pt}&&\qquad\le L_3(x,y)\| y-x\|_-\cdot \bigl( \bigl\|\nabla V(y)-\nabla V(x)\bigr\|_-+
\bigl\| \nabla U(y)+\nabla U(x) \bigr\|_- \bigr)\cdot\| \xi\|_-
\\
\hspace*{-4pt}&&\qquad\le L_3(x,y)\| y-x\|_-\cdot \bigl( L_2(x,y)\| y-x
\|_-+\bigl\|\nabla U(y)+\nabla U(x) \bigr\|_- \bigr)\cdot\| \xi\|_-,
\end{eqnarray*}
and, correspondingly,
\begin{eqnarray*}
&&\bigl|\bigl(\nabla V(y)-\nabla V(x)+\nabla U(y)+\nabla U(x)\bigr)\cdot
\nabla ^2V(y)\cdot \bigl(I+\nabla^2V(x)\bigr)\cdot\xi\bigr|
\\
&&\qquad\le L_2(y,y)\cdot \bigl(L_2(x,y)\| y-x\|_-+\bigl\|\nabla
U(y)+\nabla U(x) \bigr\|_- \bigr)\\
&&\qquad\quad{}\times\bigl(1+L_2(x,x)\bigr)\| \xi\|_-.
\end{eqnarray*}
\upqed\end{pf}
%
%We also recall the estimate
%from Lemma~\ref{lemG} (3).
By combining Proposition~\ref{propK} with the estimates in Lemma~\ref
{lemL} and Lemma~\ref{lemH}, we will now prove Proposition~\ref{propM}.

%If Assumption~\ref{assA1} is satisfied then there exists a polynomial
%$Q:\mathbb R^2\to\mathbb R$ of degree $\max(p_4+3,p_3+p_2+2,3p_2+1)$
%such that
%
%Again the polynomial $Q$ is explicit. It depends only on the values
%$c_2,c_3,c_4,p_2,p_3,p_4$ in Assumption~\ref{assA1} and on the moments
%$m_n=E[\|Z\|_-^n]$ for $n\leq\max(p_4+3,p_3+p_2+2,2p_2+1)$, but it
%does not depend on the dimension $d$.

\begin{pf*}{Proof of Proposition~\ref{propM}}
Fix $h\in(0,2)$. By (\ref{eq9}) and (\ref{eq9a}), for any $x,\tilde
x\in\mathbb R^d$,
%
%
%e5.12 #&#
\begin{eqnarray}\label{eqM1}
&&\bigl\llVert \alpha_h\bigl( x, Y_h(x)\bigr)-
\alpha_h\bigl( \tilde x, Y_h(\tilde x)\bigr)\bigr\rrVert
_{L^k}
\nonumber
\\
&&\qquad\le \bigl\llVert G_h\bigl( x, Y_h(x)\bigr)-
G_h\bigl( \tilde x, Y_h(\tilde x)\bigr)\bigr\rrVert
_{L^k}
\\
&&\qquad\le\llVert x-\tilde x\rrVert _-\cdot\sup_{z\in[x,\tilde
x]}\bigl\llVert
\bigl\|\nabla_xG_h\bigl(x,Y_h(x)\bigr)\bigr\|_+ \bigr
\rrVert _{L^k}.\nonumber
\end{eqnarray}
Moreover, by (\ref{eq45}) and Proposition~\ref{propK},
%
%
%e5.13 #&#
\begin{eqnarray}\label{eq53}
&&\bigl\llVert \bigl\|\nabla_xG_h\bigl(x,Y_h(x)
\bigr)\bigr\|_+ \bigr\rrVert _{L^k}\nonumber\\
&&\qquad=\bigl\| \bigl\|\nabla _xF_h
\bigl(x,\sqrt{h-h^2/4} Z\bigr) \bigr\|_+ \bigr\|_{L^k}
\\
&&\qquad\leq I+\frac h4\cdot \mathit{II}+\frac{h}{8-2h}\cdot \mathit{III}+\frac
{h^2}{16-4h}\cdot \mathit{IV},\nonumber
\end{eqnarray}
where
\begin{eqnarray*}
I&=&\mathbb{E}\bigl[\bigl\|\nabla V\bigl(Y_h(x)\bigr)-\nabla V(x)-
\tfrac12\bigl(Y_h(x)-x\bigr)\\
&&\hspace*{59pt}{}\times \bigl(\nabla^2V
\bigl(Y_h(x)\bigr)+\nabla^2V(x)\bigr)\bigr\|_+^k
\bigr]^{1/k},
\\
\mathit{II}&=&\mathbb{E}\bigl[\bigl\|\bigl(Y_h(x)-x\bigr)\cdot
\nabla^2V\bigl(Y_h(x)\bigr)\cdot\bigl(I+\nabla
^2V(x)\bigr)\bigr\|_+^k\bigr]^{1/k},
\\
\mathit{III}&=&\mathbb{E}\bigl[\bigl\|\bigl(\nabla V\bigl(Y_h(x)\bigr)-\nabla
V(x)+\nabla U\bigl(Y_h(x)\bigr)+\nabla U(x)\bigr)
\\
&&\hspace*{102pt}{} \times\bigl(\nabla^2V\bigl(Y_h(x)\bigr)-
\nabla^2 V(x)\bigr)\bigr\|_+^k\bigr]^{1/k},
\\
\mathit{IV}&=&\mathbb{E}\bigl[\bigl\|\bigl(\nabla V\bigl(Y_h(x)\bigr)-\nabla
V(x)+\nabla U\bigl(Y_h(x)\bigr)+\nabla U(x)\bigr)
\\
&&\hspace*{88pt}{} \times\nabla^2V\bigl(Y_h(x)\bigr)\cdot\bigl(I+
\nabla^2V(x)\bigr)\bigr\|_+^k\bigr]^{1/k}.
\end{eqnarray*}
By applying the estimates from Lemmas \ref{lemL} and~\ref{lemG}(4), we obtain
%
%
%e5.14 #&#
%e5.15 #&#
%e5.16 #&#
%e5.17 #&#
\begin{eqnarray}
\label{eq54} I&\leq&\tfrac{1}{12}\mathbb{E} \bigl[L_4
\bigl(x,Y_h(x)\bigr)^k\bigl\| Y_h(x)-x
\bigr\|_-^{3k} \bigr]^{1/k},
\\
\label{eq55} \mathit{II}&\leq&\bigl(1+L_2(x,x)\bigr)\cdot\mathbb{E}
\bigl[L_2\bigl(Y_h(x),Y_h(x)
\bigr)^k\bigl\|Y_h(x)-x\bigr\|_-^k
\bigr]^{1/k},
\\
\qquad\mathit{III}&\leq&\mathbb{E} \bigl[L_3\bigl(x,Y_h(x)
\bigr)^k\bigl\|Y_h(x)-x\bigr\|_-^k
\nonumber
\\[-8pt]
\\[-8pt]
\nonumber
\label{eq56} &&\hspace*{9pt}{} \times\bigl(\bigl(1+2L_2\bigl(x,Y_h(x)
\bigr)\bigr)^k\bigl\|Y_h(x)-x\bigr\|_-^k+2\bigl\| \nabla
U(x)\bigr\|_-^k\bigr) \bigr]^{1/k},
\\
\label{eq57} \mathit{IV}&\leq&\bigl(1+L_2(x,x)\bigr)
\nonumber
\\
&&{}\times\mathbb{E}
\bigl[L_2\bigl(Y_h(x),Y_h(x)
\bigr)^k\\
&&\hspace*{22pt}{}\times\bigl(\bigl(1+L_2\bigl(x,Y_h(x)
\bigr)\bigr)^k\bigl\|Y_h(x)-x\bigr\|_-^k+2\bigl\| \nabla
U(x)\bigr\|_-^k\bigr) \bigr]^{1/k}.\nonumber
\end{eqnarray}
The assertion for semi-implicit {E}uler proposals is now a direct
consequence of the estimates in Lemma~\ref{lemH}, (\ref{eq40}) and
(\ref
{eqM1}).

The assertion for Ornstein--Uhlenbeck
proposals follows in a similar way from~(\ref{eqKOU}), Lemma~\ref{lemL}(1) and
the estimates in Lemma~\ref{lemH}.
\end{pf*}

Again, it is possible to write down the polynomial in Proposition~\ref
{propM} explicitly.
For semi-implicit {E}uler proposals, we illustrate this in the case
$k=1$ and $p_2=p_3=p_4=0$. Here, by (\ref{eq54}), (\ref{eq55}), (\ref
{eq56}) and (\ref{eq57}) we obtain
\begin{eqnarray*}
I&\leq&\frac{C_4}{12}\mathbb{E} \biggl[ \biggl(\frac h2\bigl\|\nabla U(x)\bigr\| _-+
\sqrt h\|Z\|_- \biggr)^3 \biggr] \leq \frac{C_4}{4} \biggl(
\frac{h^3}{8}\bigl\|\nabla U(x)\bigr\| ^3_-+h^{3/2}m_3
\biggr),
\\
\mathit{II}&\leq&\bigl(C_2+C_2^2\bigr) \mathbb{E}
\biggl[\frac{h}{2}\bigl\|\nabla U(x)\bigr\| _-+\sqrt h\|Z\|_- \biggr]
\\
&=& \bigl(C_2+C_2^2\bigr) \biggl(\frac h2\bigl\|
\nabla U(x)\bigr\|_-+h^{1/2}m_1 \biggr),
\\
\mathit{III}&\leq&C_3 \biggl(2\bigl\|\nabla U(x)\bigr\|_- \cdot\mathbb{E} \biggl[\frac
h2\bigl\| \nabla U(x)\bigr\|_-+\sqrt h\|Z\|_- \biggr]
\\
&&\hspace*{18pt}{} +(1+2C_2) \mathbb{E} \biggl[ \biggl(\frac h2\bigl\|\nabla U(x)
\bigr\| _-+\sqrt h\|Z\|_- \biggr)^2 \biggr] \biggr)
\\
&\leq&2C_3\bigl\|\nabla U(x)\bigr\|_- \biggl(\frac h2 \bigl\|\nabla U(x)\bigr\|_-+\sqrt
h m_1 \biggr)
\\
&&{} +C_3(1+2C_2) \biggl(\frac{h^2}2\bigl\|\nabla U(x)
\bigr\|^2+2hm_2 \biggr),\\
\mathit{IV}&\leq&(1+C_2)C_2 \biggl((1+C_2)\mathbb{E}
\biggl[\frac h2\bigl\|\nabla U(x)\bigr\| _-+\sqrt h\|Z\|_- \biggr]+2\bigl\|\nabla U(x)\bigr\|_-
\biggr)
\\
&\leq&2(1+C_2)C_2\bigl\|\nabla U(x)\bigr\|_-+(1+C_2)^2C_2
\biggl(\frac h2\bigl\|\nabla U(x)\bigr\|_-+\sqrt hm_1 \biggr).
\end{eqnarray*}
Hence by (\ref{eq53}), for $h\in(0,2)$,
%
%
%e5.18 #&#
\begin{eqnarray}\label{eqQE}
&&\mathbb{E}\bigl[\bigl\|\nabla_xG_h
\bigl(x,Y_h(x)\bigr)\bigr\|_+\bigr] \nonumber\\
&&\qquad\leq\tfrac14 h^{3/2}
\bigl(C_4m_3+(1+C_2)C_2m_1+2C_3
\bigl\|\nabla U(x)\bigr\|_-m_1\bigr)
\nonumber\\
&&\qquad\quad{}+\tfrac18h^2\bigl(4C_2(1+2C_2)m_2+3C_2(1+C_2)
\bigl\|\nabla U(x)\bigr\|_-+2C_3\bigl\| \nabla U(x)\bigr\|_-^2\bigr)
\\
&&\qquad\quad{}+\tfrac1{16}h^{5/2}C_2(1+C_2)^2
\bigl(2m_1+h^{1/2}\bigl\|\nabla U(x)\bigr\|_-\bigr)
\nonumber\\
&&\qquad\quad{}+\tfrac{1}{32}h^3 \bigl(4C_3(1+2C_2)
\bigl\|\nabla U(x)\bigr\|_-^2+C_4\bigl\|\nabla U(x)\bigr\|_-^3
\bigr).\nonumber
\end{eqnarray}
%
%By Lemma~\ref{lemA} and Proposition~\ref{propM} we obtain the
%following result on the dependence of the acceptance probabilities on
%the current state:

For \emph{Ornstein--Uhlenbeck proposals}, we prove the explicit bound
for the
dependence of the rejection probabilities on the current state for the case
$k=2$ and $p_2=0$ as stated in Proposition~\ref{propMOU}.
\begin{pf*}{Proof of Proposition~\ref{propMOU}, second part}
If $p_2=0$, then by (\ref{eq45OU}), (\ref{eqKOU}) and (\ref{eqnablav}),
\begin{eqnarray*}
\bigl\|\nabla_xG^\mathrm{OU}\bigl(x,Y_h^\mathrm{OU}(x)
\bigr)\bigr\|_+ &\le& \bigl\| \nabla V\bigl(Y_h^\mathrm{OU}(x)\bigr)-
\nabla V(x)\bigr\|_+ + \frac h2\bigl \| \nabla V(x)\bigr\|_+
\\
&\le& C_2\bigl\| Y_h^\mathrm{OU}(x)-x\bigr\|_- +
\bigl(C_1+C_2\| x\|_-\bigr) h/2
\\
& \le& C_2\| Z\|_- h^{1/2} + \bigl(C_1+2C_2
\| x\|_-\bigr) h/2
\end{eqnarray*}
for any $x\in\mathbb R^d$. Therefore,
\[
\mathbb{E} \bigl[\bigl\|\nabla_xG^\mathrm{OU}
\bigl(x,Y_h^\mathrm{OU}(x)\bigr)\bigr\| _+^2
\bigr]^{1/2} \le C_2m_2^{1/2}
h^{1/2} + \bigl(C_1+2C_2\| x\|_-\bigr) h/2.
\]
The assertion now follows similarly to (\ref{eqM1}).
\end{pf*}

%s6 #&#
\section{Upper bound for exit probabilities}
\label{secLyapunov}
In this section, we prove an upper bound for the exit probabilities of
the MALA chain
from the ball $B_R^-$ that is required in the proof of Theorem~\ref
{thmMAIN}; cf. \cite{G} for a detailed proof of a more general result. Let
%
%
%e6.1 #&#
\begin{equation}
\label{eq6A} f(x) := \exp \bigl(K\| x\|_-^2/16 \bigr).
\end{equation}
The following lemma shows that $f(x)$ acts as a Lyapunov function for
the MALA
transition kernel on $B_R^-$:

%
%le6.1 #&#
\begin{lemma}\label{lem6a}
Suppose that Assumptions \ref{assA1} and \ref{assA2} hold. Then there
exist constants
$C_1,C_2,\rho_1\in(0,\infty)$ such that
%
%
%e6.2 #&#
\begin{equation}
\label{eq6B} q_hf \le f^{1-Kh/4} e^{C_2h} \qquad\mbox{on
}B_R^-
\end{equation}
for any $R,h\in(0,\infty)$ such that $h^{-1}\ge C_1(1+R)^{\rho_1}$.
\end{lemma}

\begin{pf}
We first observe that a corresponding bound holds for the proposal
kernel $p_h$. Indeed,
by (\ref{eq32}), and since $\| v\|_-^2=v\cdot Gv$ with a nonnegative
definite symmetric
matrix $G\le I$, an explicit computation yields
\begin{eqnarray*}
(p_hf) (x) &=& \E \biggl[ \exp\biggl(K\biggl\| x-\frac h2\nabla U(x)+
\sqrt {h-h^2/4}Z\biggr\|_-^2\Big/16\biggr) \biggr]
\\
&\le& \exp \biggl( K(1+Kh/4)\biggl\| x-\frac h2\nabla U(x)\biggr\|_-^2\Big/16
\biggr).
\end{eqnarray*}
Moreover, by Assumption~\ref{assA2},
\[
\biggl\| x-\frac h2\nabla U(x)\biggr\|_-^2 \le \biggl(1-\frac{Kh}2
\biggr)\| x\|_-^2+\frac h{2K}\bigl\|\nabla U(0)\bigr\|_-^2+
\frac{h^2}4\bigl\|\nabla U(x)\bigr\|_-^2.
\]
Hence by Assumption~\ref{assA1},
there exist constants $C_3,C_4,\rho_2\in
(0,\infty)$ such that
\[
(p_hf) (x) \le f(x)^{1-Kh/4} e^{C_3h}
\]
for any $x\in\mathbb R^d$ and $h\in(0,\infty)$ such that $h^{-1}\ge
C_4 (1+\| x\|_-)^{\rho_2}$. By the upper bound for the rejection
probabilities in
Proposition~\ref{propI}, we conclude that there exists a polynomial $s$
such that the
corresponding upper bound
\begin{eqnarray*}
q_hf &\le& f^{1-Kh/4}e^{C_3h} + s(R)h^{3/2}
f \\
&= &f^{1-Kh/4} \bigl(e^{C_3h} + s(R)h^{3/2}
f^{Kh/4} \bigr)
\\
&\le& f^{1-Kh/4}e^{(C_3+1)h}
\end{eqnarray*}
holds on $B_R^-$ whenever both $h^{-1}\ge C_4 (1+R)^{\rho_2}$ and
$s(R)h^{1/2}f^{Kh/4}\le1 $. The assertion follows, since the second
condition is
satisfied if $K^2hR^2/64\le1$ and $s(R)eh^{1/2}\le1$.
\end{pf}

Now consider the first exit time
\[
T_R := \inf\bigl\{ n\ge0\dvtx X_n\notin B_R^-
\bigr\},
\]
where $(X_n,\P_x)$ is the Markov chain with transition kernel $q_h$ and
initial condition
$X_0=x$ $P_x$-a.s. We can estimate $T_R$ by constructing a
supermartingale based
on the Lyapunov function $f$:

%
%th6.2 #&#
\begin{thmm}\label{thmm6b}
If Assumptions \ref{assA1} and \ref{assA2} hold, then there exist
constants $C,\rho,
D\in(0,\infty)$ such that the upper bound
%
%
%e6.3 #&#
\begin{equation}
\label{eq6c} \P_x [T_R\le n] \le D nh \exp\bigl[K\bigl(
\| x\|_-^2-R^2\bigr)/24\bigr]
\end{equation}
holds for any $n \ge0$, $R,h\in(0,\infty)$ such that $h^{-1}\ge
C(1+R)^\rho$,
and $x\in B_R^-$.
\end{thmm}

\begin{pf}
Fix $n\in\mathbb N $, choose $C_1,C_2,\rho_1$ as in the lemma above,
and let
\[
M_j := f(X_j)^{(1-Kh/4)^{n-j}} \exp \Biggl(
-C_2h\sum_{i=0}^j(1-Kh/4)^{n-i}
\Biggr)
\]
for $j=0,1,\ldots,n$. If $h^{-1}\ge C_1(1+R)^{\rho_1}$, then by
Jensen's inequality and
(\ref{eq6B}),
\begin{eqnarray*}
\E_x [ M_{j+1}|\mathcal F_j ] &\le&
(q_hf) (X_j)^{1-Kh/4}\exp \Biggl(
-C_2h\sum_{i=0}^{j+1}(1-Kh/4)^{n-i}
\Biggr)
\\
&\le& M_j \qquad\mbox{on } \bigl\{ X_j\in B_R^-
\bigr\} \mbox{ for any }j<n.
\end{eqnarray*}
Hence the stopped process $(M_j^{T_R})_{0\le j\le n}$ is a
supermartingale, and thus
\[
\E_x [ M_{T_R}; n-m\le T_R\le n ] \le
\E_x [ M_0 ] \qquad\mbox{for any }0\le m\le n.
\]
Noting that $M_0=f(x)^{(1-Kh/4)^n}=\exp ( (1-Kh/4)^nK\| x\|
_-^2/16 ) $, and
\begin{eqnarray*}
M_{T_R} &\ge& \bigl(f(X_{T_R})\exp(-4C_2/K)
\bigr)^{(1-Kh/4)^{n-T_R}}
\\
&= & \exp \biggl[ \biggl(\frac{K}{16}R^2-4C_2/K
\biggr)\cdot(1-Kh/4)^{n-T_R} \biggr],
\end{eqnarray*}
we obtain the bound
\begin{eqnarray*}
&&\P_x[n-m\le T_R\le n] \\
&&\qquad\le \exp \biggl[ \biggl(1-
\frac{Kh}4 \biggr)^m \biggl(\frac{K}{16}\biggl(\biggl(1-
\frac{Kh}4\biggr)^{n-m}\| x\|_-^2-R^2
\biggr)+\frac
{4C_2}K \biggr) \biggr]
\end{eqnarray*}
for any $0\le m\le n$ provided $R^2\ge64 C_2/K^2$. In particular, if
$mKh/2\le
\log(3/2)$, then $(1-Kh/4)^m\ge\exp(-mKh/2)\ge2/3$, and hence
\[
\P_x [n-m\le T_R\le n] \le \exp(4 C_2K)
\cdot\exp \bigl( K\bigl(\| x\| _-^2-R^2\bigr)/24 \bigr).
\]
The assertion follows by partitioning $\{ 0,1,\ldots,n\} $ into blocks
of length $\le m$
where $m=\lfloor2\log(3/2) K^{-1}h^{-1}\rfloor$.
\end{pf}

%s7 #&#
\section{Proof of the main results}
\label{secFinal} In this section, we combine the results in order to
derive the contraction properties for Metropolis--Hastings transition
kernels stated in Theorems \ref{thmQ}, \ref{thmQOU} and \ref{thmMAIN},
and we finally prove \ref{thmFINAL}. Note that for $x,\tilde
x\in\mathbb{R}^d$, the distances
%
%
%e7.1 #&#
%e7.2 #&#
\begin{eqnarray}
\bigl\llVert Y_h^{\mathrm{OU}}(x)-Y_h^{\mathrm{OU}}(
\tilde x)\bigr\rrVert _- &=& (1-h/2 ) \|x-\tilde x\|_- \qquad\mbox{and}\label{eqdistOU}
\\
\bigl\llVert Y_h(x)-Y_h(\tilde x)\bigr\rrVert _- &=&
\bigl\| x-\tilde x-\bigl(\nabla U(x)-\nabla U(\tilde x)\bigr) h/2\bigr\|_-
\end{eqnarray}
are deterministic. We now combine Lemma~\ref{lemD} with the estimates in
Propositions~\ref{propI} and \ref{propM}:

\begin{pf*}{Proof of Theorem~\ref{thmQ}}
We fix $h\in(0,2)$, $R\in
(0,\infty)$
and $x,\tilde x\in B_R^-$. By the basic contractivity Lemma \ref{lemD}
and by (\ref{eqlemD}), respectively,
\begin{eqnarray*}
&&\mathbb{E} \bigl[ \bigl\llVert W_h(x)-W_h(
\tilde x)\bigr\rrVert _- \bigr]\\
&&\qquad\le \| x-\tilde x\|_-
\\
&&\qquad\quad{}- \bigl(1-\E \bigl[\max\bigl(1-\alpha_h \bigl(x,Y_h(x)
\bigr),1-\alpha_h \bigl(\tilde x, Y_h(\tilde x)\bigr)
\bigr) \bigr] \bigr)\\
&&\qquad\quad{}\times  \bigl(\|x-\tilde x\|_- -\bigl\|Y_h(x)-Y_h(
\tilde x)\bigr\|_- \bigr)
\\
&&\qquad\quad{}+ \mathbb{E}\bigl[\max\bigl(\bigl\|x-Y_h(x)\bigr\|_-, \bigl\|\tilde
x-Y_h(\tilde x)\bigr\|_-\bigr)^2\bigr]^{1/2} \\
&&\qquad\quad{}\times
\mathbb{E}\bigl[\bigl(\alpha_h \bigl(x,Y_h(x)\bigr)-
\alpha_h \bigl(\tilde x, Y_h(\tilde x)\bigr)
\bigr)^2\bigr]^{1/2}.
\end{eqnarray*}
By Proposition~\ref{propP},
\[
\bigl\|Y_h(x)-Y_h(\tilde x)\bigr\|_- \le \bigl(1-Kh/2+M(R)^2h^2/8
\bigr)\cdot\| x-\tilde x\|_-,
\]
and by Lemma~\ref{lemH} (2),
\begin{eqnarray*}
&&\mathbb{E}\bigl[\max\bigl(\bigl\|x-Y_h(x)\bigr\|_-, \bigl\|\tilde
x-Y_h(\tilde x)\bigr\|_-\bigr)^2\bigr]^{1/2} \\
&&\qquad\le
m_2^{1/2} h^{1/2}+ \max\bigl(\bigl\|\nabla U(x)\bigr\|_-, \bigl\|
\nabla U(\tilde x)\bigr\|_-\bigr) h/2.
\end{eqnarray*}
The assertion of Theorem~\ref{thmQ} follows by combining these
estimates with
the bounds for the acceptance probabilities in Propositions \ref
{propI} and
\ref{propM}.
\end{pf*}

The corresponding bound for Ornstein--Uhlenbeck proposals follows
similarly from Lemma~\ref{lemD} and Proposition~\ref{propMOU}:

\begin{pf*}{Proof of Theorem~\ref{thmQOU}}
We again fix $h\in(0,2)$,
$R\in(0,\infty)$,
and $x,\tilde x\in B_R^-$. Since $Y_h^\mathrm{OU}(x)-Y_h^\mathrm
{OU}(\tilde x)=(1-h/2)(x-\tilde x )$ and $\| x-Y_h^\mathrm{OU}(x)\|
_-\le\| x\|_-h/2 +\| Z\|_-\sqrt h$, the basic contractivity Lemma \ref
{lemD} implies
\begin{eqnarray*}
&&\mathbb{E} \bigl[ \bigl\llVert W_h^\mathrm{OU}(x)-W_h^\mathrm
{OU}(\tilde x)\bigr\rrVert _- \bigr] \\
&&\qquad\le \biggl(1-\frac h2\biggr) \| x-\tilde x
\|_-
\\
&&\qquad\quad{}+ \frac h2\|x-\tilde x\|_- \E \bigl[\max\bigl(1-\alpha^\mathrm{OU}
\bigl(x,Y_h^\mathrm{OU}(x)\bigr),1-\alpha^\mathrm{OU}
\bigl(\tilde x, Y_h^\mathrm {OU}(\tilde x)\bigr)\bigr) \bigr]
\\
&&\qquad\quad{}+ \biggl( \frac h2\max\bigl(\| x\|_-,\| \tilde x\|_-\bigr)+\sqrt h m_2^{1/2}
\biggr)\\
&&\qquad\quad{}\times \mathbb{E}\bigl[\bigl(\alpha^\mathrm{OU} \bigl(x,Y_h^\mathrm{OU}(x)
\bigr)-\alpha ^\mathrm{OU} \bigl(\tilde x, Y_h^\mathrm{OU}(
\tilde x)\bigr)\bigr)^2\bigr]^{1/2}.
\end{eqnarray*}
The assertion of Theorem~\ref{thmQOU} follows by combining this
estimates with
the bounds for the acceptance probabilities in Proposition~\ref{propMOU}.
\end{pf*}

\begin{pf*}{Proof of Theorem~\ref{thmMAIN}}
Noting that
\[
\| x\|_- -(2R)^2 \le -3R^2 \qquad\mbox{for any } x\in
B_R^-,
\]
the assertion is a direct consequence of Corollary~\ref{corWbound} and
Theorem~\ref{thmm6b} applied with $R$ replaced by $2R$.
\end{pf*}

Let $\mu_R=\mu(\cdot|B_R^-)$ denote the conditional measure on $B_R^-$.
The fact that $\mu_R$ is a stationary distribution for the Metropolis--Hastings
transition kernel $q_h$ can be used to bound the Wasserstein distance between
$\mu_Rq_h^n$ and $\mu_R$:

%
%le7.1 #&#
\begin{lemma}\label{lem7a}
For any $R\ge0$ and $a\in(0,1)$,
\begin{eqnarray*}
\mathcal W_{2R}\bigl(\mu_Rq_h^n,
\mu_R\bigr) &\le& 8R \bigl(\mu_Rq_h^n
\bigr) \bigl(\mathbb R^d\setminus B_R^-\bigr)
\\
&\le& 8{(1-a)^{-1}} \mathcal W_{2R}\bigl(
\mu_Rq_h^n,\delta_0q_h^n
\bigr) + 8R \bigl(\delta_0q_h^n\bigr)
\bigl(B_{aR}^-\bigr).
\end{eqnarray*}
\end{lemma}

\begin{pf}
The distance induced by the total variation norm $\|\cdot\|_{\mathrm{TV}}$ is
the Wasserstein distance
w.r.t. the metric $d(x,y)=I_{\{ x\neq y\} }$. Since $d_{2R}(x,y)\le
4Rd(x,y)$, we obtain
%
%
%e7.3 #&#
\begin{eqnarray}\label{eq7x}
\mathcal W_{2R}\bigl(\mu_Rq_h^n,
\mu_R\bigr) &\le& 4R\bigl\|\mu_Rq_h^n-
\mu_R\bigr\|_{\mathrm{TV}}\nonumber\\
& =& 8R\bigl\|\bigl(\mu_Rq_h^n-
\mu_R\bigr)^+\bigr\|_{\mathrm{TV}}
\\
&\le&8R \bigl(\mu_Rq_h^n\bigr) \bigl(
\mathbb R^d\setminus B_R^-\bigr).\nonumber
\end{eqnarray}
Here we have used in the last step that $\mu q_h=\mu$, and hence
\begin{eqnarray*}
\bigl(\mu_Rq_h^n\bigr) (A) &\le&\bigl(
\mu_Rq_h^n\bigr) \bigl(A\cap
B_R^-\bigr)+\bigl(\mu_Rq_h^n
\bigr) \bigl(\mathbb R^d\setminus B_R^-\bigr)
\\
&\le&\mu_R(A)+\bigl(\mu_Rq_h^n
\bigr) \bigl(\mathbb R^d\setminus B_R^-\bigr)
\end{eqnarray*}
for any Borel set $A\subseteq\mathbb R^d$. Moreover, for $a\in(0,1)$,
%
%
%e7.4 #&#
\begin{eqnarray}
\label{eq7y} &&\mathcal W_{2R}\bigl(\mu_Rq_h^n,
\delta_0q_h^n\bigr)
\nonumber
\\[-8pt]
\\[-8pt]
\nonumber
&&\qquad\ge (R-aR)\cdot \bigl(
\bigl(\mu_Rq_h^n\bigr) \bigl(\mathbb
R^d\setminus B_R^-\bigr)-\bigl(\delta
_0q_h^n\bigr) \bigl(\mathbb R^d
\setminus B_{aR}^-\bigr) \bigr).
\end{eqnarray}
Indeed, for any coupling $\eta(dx \,d\tilde x )$ of the two measures,
\[
\eta\bigl(d_{2R}(x,\tilde x)\ge R-aR\bigr) \ge \bigl(
\mu_Rq_h^n\bigr) \bigl(\mathbb
R^d\setminus B_R^-\bigr)-\bigl(\delta_0q_h^n
\bigr) \bigl(\mathbb R^d\setminus B_{aR}^-\bigr).
\]
The assertion follows by combining the estimates in (\ref{eq7x}) and
(\ref{eq7y}).
\end{pf}

\begin{pf*}{Proof of Theorem~\ref{thmFINAL}}
By combining the estimates in Theorem~\ref{thmMAIN}, Lemma~\ref{lem7a}
with $a=6/7$, and Theorem~\ref{thmm6b}, we obtain
\begin{eqnarray*}
\mathcal W_{2R}\bigl(\nu q_h^n,
\mu_R\bigr) &\le& \mathcal W_{2R}\bigl(\nu
q_h^n,\mu_Rq_h^n
\bigr) +\mathcal W_{2R}\bigl(\mu_R q_h^n,
\mu_R\bigr)
\\
&\le& \biggl(1-\frac K4 h\biggr)^n\mathcal W_{2R}(\nu,
\mu_R) + DR\exp\bigl(-KR^2/8\bigr) nh
\\
&&{}+56\cdot\biggl(1-\frac K4 h\biggr)^n\mathcal W_{2R}(
\mu_R,\delta_0) + 56 DR\exp\bigl(-KR^2/8
\bigr) nh\\
&&{}+8R \mathbb P_0[T_{6R/7}\le n]
\\
&\le&58 R\cdot\biggl(1-\frac K4 h\biggr)^n + 57 DR\exp
\bigl(-KR^2/8\bigr) nh\\
&&{}+8 DR\exp\bigl(-KR^2/33\bigr) nh.
\end{eqnarray*}
\upqed\end{pf*}

\section*{Acknowledgments}
I would like to thank in particular
Daniel Gruhlke and Sebastian Vollmer for many fruitful discussions related
to the contents of this paper.

%%\bibliography{mala}
%% \MRhref is called by the amsart/book/proc definition of \MR.
% \href{http://www.ams.org/mathscinet-getitem?mr=#1}{#2}
%}
%
% imsref loaded by akundreckaite, 2013-09-04 12:53:55
% imsref loaded by akundreckaite, 2013-09-04 13:40:35
%


%suskaldyti doi

\printaddresses

\end{document}